\documentclass[11pt,a4paper]{article}
\usepackage{amsfonts}
\usepackage{amsmath,amssymb}
\usepackage{mathrsfs}
\usepackage{hyperref}
\usepackage{graphicx}
\usepackage{indentfirst}
\usepackage{color}

\newtheorem{thm}{Theorem}[section]

\newtheorem{lem}[thm]{Lemma}
\newtheorem{prop}[thm]{Proposition}

\newtheorem{rem}[thm]{Remark}

\numberwithin{equation}{section}\allowdisplaybreaks

\def\leq{\leqslant}

\def\leq{\leqslant}
\def\geq{\geqslant}

\newcommand\blfootnote[1]{%
 \begingroup
 \renewcommand\thefootnote{}\footnote{#1}%
 \addtocounter{footnote}{-1}%
 \endgroup
 }

\begin{document}

\title{\Large\bf  On Dissipative Nonlinear Evolutional Pseudo-Differential Equations}

\author{\normalsize   Mingjuan Chen, \ \   Baoxiang Wang, \ \ Shuxia Wang,  \ \ M. W. Wong\blfootnote{Baoxiang Wang is the corresponding author.  Mingjuan Chen and   Baoxiang Wang, LMAM, School of Mathematical Sciences, Peking University, Beijing 100871, China,
Emails: mjchenhappy@pku.edu.cn, \  wbx@math.pku.edu.cn;  Shuxia Wang,  Department of Mathematics, Beijing Institute of Petrochemical Technology, Beijing 102617, China, Email: wangshuxia@bipt.edu.cn;  M. W. Wong, Department of Mathematics and Statistics,
 York University,
4700 Keele Street,
Toronto,
Canada,  Email: mwwong@mathstat.yorku.ca.
} \\
} \maketitle

\thispagestyle{empty}
\begin{abstract}\footnotesize

First, using the uniform decomposition in both physical and frequency spaces, we obtain an equivalent norm on modulation spaces. Secondly, we consider the Cauchy problem for the dissipative evolutionary pseudo-differential equation
\begin{align}
\partial_t u  + A(x,D)  u = F\big((\partial^\alpha_x u)_{|\alpha|\leq \kappa}\big), \ \     u(0,x)=   u_0(x), \nonumber
\end{align}
where $A(x,D)$ is a dissipative pseudo-differential operator  and $F(z)$ is a multi-polynomial.  We will develop the uniform decomposition techniques in both physical and frequency spaces to study its local well posedness in modulation spaces $M^s_{p,q}$ and in Sobolev spaces $H^s$. Moreover, the local solution can be extended to a global one in $L^2$ and in $H^s$ ($s>\kappa+d/2$) for certain nonlinearities.\\

{\bf MSC 2010:} 35S30,  42B37, 42B35, 35K55.

\end{abstract}

\section{Introduction}

We study the Cauchy problem of  the dissipative evolutionary pseudo-differential equation
\begin{align}
\partial_t u  + A(x,D)  u = F\big((\partial^\alpha_x u)_{|\alpha|\leq \kappa}\big), \ \     u(0,x)=   u_0, \label{PDE}
\end{align}
where $A(x,D)$ is a pseudo-differential operator, which is defined by (cf. \cite{St93, Wo14})
\begin{align*}
A(x,D)f(x)=\int e^{{\rm i}x\cdot\xi}A(x,\xi)\widehat{f}(\xi)d\xi, \ \ \ \ f\in\mathscr{S}.
\end{align*}
Recall that
$A(\cdot,\ \cdot)$ is said to be in the class $S^M$ if $A \in C^N(\mathbb{R}^{2d})$ for some $N\gg d$ and satisfies
\begin{align*}
|\partial_x^\alpha\partial_\xi^\beta A(x,\xi)|\leq A_{\alpha,\beta}(1+|\xi|)^{M-|\beta|}, \ \ x,\xi \in \mathbb{R}^d, \ \ |\alpha|, |\beta| \leq N.
\end{align*}
Moreover, $A(x,D)$ is said to be dissipative with order $M$, if
\begin{align*}
 \mathfrak{Re}\ A(x,\xi) \gtrsim |\xi|^M, \ \ (x,\xi)\in\mathbb{R}^d\times\mathbb{R}^d.
\end{align*}

We will study the following two cases: (i) $A(x,D) = a(x) +b(D)$, where $\mathfrak{Re} a(x) \gtrsim | x|^{\sigma_1}$ and $b\in S^{M}$ is $M$-order dissipative for some $\sigma_1,M>0$,  (ii)  $A\in S^M$  is $M$-order dissipative for some $M>0$. The nonlinearity
\begin{align} \label{nonlinearity}
F(z) = \sum_{2\leq |\beta| \leq K} \lambda_\beta z^\beta, \ \  z=((\partial^\alpha_x u)_{|\alpha|\leq \kappa})
\end{align}
is a multi-polynomial which contains at most $\kappa$-order derivative of $u$, $\kappa $ is less than or equal to the dissipative order $M$, the powers in the nonlinearity are less than or equal to $K \in \mathbb{N}$.

In order to solve the equation \eqref{PDE}, we will use   the uniform decomposition operators to both frequency and physical spaces. The uniform decomposition in physical spaces was applied by Kenig, Ponce and Vega \cite{KPV93, KPV98} in the study of the derivative nonlinear Schr\"odinger equations. The frequency-uniform decomposition operators have been extensively applied in the study of nonlinear evolution equations (cf.
\cite{Wa07,Wa07a,Wa06}). First, we recall the definition of frequency-uniform decomposition operators. Let $\{\sigma_n\}_{n\in \mathbb{Z}^d}$ be
a smoothing function sequence satisfying
$$\left\{
\begin{array}{l}
   {\rm supp} \ \sigma \subset [-3/4,3/4]^d, \ \sigma_n(\cdot)= \sigma (\cdot-n),\\
  \sum_{n\in \mathbb{Z}^d} \sigma_n (\xi) \equiv 1, \quad \forall \;
\xi \in
\mathbb{R}^d, \\
  |\partial^{\alpha} \sigma_n (\xi)|   \leq C^{|\alpha|}, \ \  \forall \;
\xi \in
\mathbb{R}^d.
\end{array}
\right. \eqno{(\rm UD)}
$$
Denote
\begin{align*}
\Upsilon_d = \left\{ \{\sigma_n\}_{n\in \mathbb{Z}^d}: \;
\{\sigma_n\}_{n\in \mathbb{Z}^d}\;\; \mbox {satisfies (UD)}
\right\}.
\end{align*}
Let
$\{\sigma_n\}_{n\in \mathbb{Z}^d}\in \Upsilon_d$ be a function
sequence and
\begin{align*}
\Box_n := \mathscr{F}^{-1} \sigma_n(\xi) \mathscr{F},
\end{align*}
$\Box_n$ ($n\in \mathbb{
Z}^d$) are said to be the frequency-uniform decomposition operators. For $n\in \mathbb{Z}^d$, we denote $\langle n\rangle=(1+|n|^2)^{1/2}$. Let $s\in \mathbb{R}$, $1\leq p,q\leq \infty$. Modulation spaces $M^s_{p,q}$ were introduced by Feichtinger \cite{Fe83} in 1983. Using the frequency-uniform decomposition operators, they can be equivalently defined as
\begin{align}
M_{p,q}^s(\mathbb{R}^d)=\Bigg\{f\in \mathscr{S}'(\mathbb{R}^d): \|f\|_{M_{p,q}^s}=\bigg(\sum_{ n\in\mathbb{Z}^d}\langle n\rangle^{sq}\|\Box_{ n}f\|_p^q\bigg)^{1/q} <\infty \Bigg\}. \label{Modulation}
\end{align}
It is well-known that $M^s_{2,2}$ is identical with Sobolev spaces $H^s:= (I-\Delta)^{-s/2 }L^2(\mathbb{R}^d)$.  Modulation spaces are basic tools in the theory of time-frequency, which plays an important role in the Gabor frame theory, see Gr\"ochenig \cite{Gr01}. At first glance, one sees that they are quite similar to Besov spaces, which can be defined by using the uniform decomposition instead of the dyadic decomposition in Besov spaces, however, comparing with Besov spaces, modulation spaces exhibit rather different aspects in their regularity, scaling properties  and algebraic structures, see \cite{Fe83, ST07, HW14}. On the other hand, modulation spaces have been widely applied in the theory of pseudo-differential operators and nonlinear evolution equations, cf. \cite{BG07, BO09, Co13, Co14, CN08, CN09, Gr01, Gr06, Gr08, KaKoIt14, Kato14, Iw10, Ni14, ST08, To04, To04b, Wa07, Wa07a, Wa06}. The solutions of the dispersive wave equations in modulation spaces also have different behavior with respect to Besov spaces. For example, $e^{it\Delta}$ is bounded in all modulation spaces and $e^{it\Delta}$ is bounded in Besov spaces spaces $B^s_{p,q}$ if and only if $p=2$.    For the nonlinear case, we know that nonlinear Schr\"odinger equation (NLS) is globally well-posed in $M^0_{2,1}$ with small data, which contains a class of functions out of the critical Sobolev spaces $H^{s_c}$, so called super-critical data in $H^s$  ($0<s<s_c$) for NLS.  The global well-posedness for such a kind of super-critical cases of NLS is hard to obtain by only using Sobolev spaces.

In order to solve equation \eqref{PDE}, we further decompose physical spaces in the same way as the frequency space and introduce the frequency-physical-uniform decomposition operators
\begin{align*}
\Box_{m,n} := \sigma_m(x) \mathscr{F}^{-1} \sigma_n(\xi) \mathscr{F}.
\end{align*}
Acted on the operators $\Box_{m,n}$, the equation \eqref{PDE} is localized in a neighbourhood of $(m,n)$ in physical and frequency spaces and one can roughly regard $A(m,n)$ as the main part of $A(x,D)$. Along this line, the dissipative structure of equation \eqref{PDE} can be simplified as $\partial_t + A(m,n) $ in local physical-frequency spaces. Indeed, we have

\begin{prop}\label{discrit} {\rm (Localization in physical and frequency spaces)}
 \eqref{PDE} is formally equivalent to the following system
\begin{align}
\label{localization}
&\partial_t \Box_{m,n} u  + A(m,n) \Box_{m,n} u  + \Box_{m,n} \mathscr{F}^{-1}(A(x,\xi)-A(m,n))\widehat{u} = \Box_{m,n}F(u),\nonumber\\
 &\ \ \ \ \ \ \ \ \    {\rm for \ all\ } m,n \in \mathbb{Z}^d,
\end{align}
where $F(u):= F\big((\partial^\alpha_x u)_{|\alpha|\leq \kappa}\big)$ (and we always use this notation if there is no confusion).
\end{prop}
{\bf Proof.} It is an easy consequence of \eqref{PDE} and (UD). $\hfill\Box$

\begin{prop}\label{discrit integral} {\rm (Localized integral equation)} Denote
\begin{align}
\label{localizationop}
\mathscr{A}_{m,n} f = \int^t_0 e^{-(t-s) A(m,n)}\Box_{m,n} f(s) ds.
\end{align}
Then \eqref{localization}  is formally equivalent to
\begin{align}
\label{localization integral}
  \Box_{m,n} u (t) = & e^{-t A(m,n)}   \Box_{m,n} u_0   - \mathscr{A}_{m,n} ( \mathscr{F}^{-1}(A(x,\xi)-A(m,n))\widehat{u})  \nonumber\\
  & +\mathscr{A}_{m,n} F(u),   \ \ {\rm for\ all\ } m,n \in \mathbb{Z}^d.
\end{align}
\end{prop}
{\bf Proof.} Considering the following equation for $c>0$,
\begin{align}\label{localization integral1}
\partial_t v  + c  v = f, \ \     v(0,x)=   v_0,
\end{align}
one easily sees that \eqref{localization integral1} is equivalent to

\begin{align*}
  v (t) = e^{-c t  }  v_0   +    \int^t_0  e^{-c(t-s)} f(s)  ds.
\end{align*}
This implies \eqref{localization integral}, as desired.  $\hfill\Box$\\

So, in view of Proposition \ref{discrit integral}, in order to show the well-posedness of \eqref{PDE}, it suffices to consider the well-posedness of the system \eqref{localization integral}. Recall that the analytical semigroup $e^{t A(x,D) }$  and its applications to semilinear heat equations in Banach (Hilbert) spaces have been studied for instance, in \cite{IW96, IW98},  but the dissipative effects of $\partial_t +A(x,D)$ were not extensively applied in those works. Noticing that $e^{-t A(m,n)}$ is exponentially decaying, the dissipative structure in $\partial_t + A(x,D)$ will be captured by the semi-group system $\{ e^{-t A(m,n)} \}_{m,n\in \mathbb{Z}^d}$  and the integral operator system $\{\mathscr{A}_{m,n} \}_{m,n\in \mathbb{Z}^d}$. Now the  difficulty is to show that $\mathscr{F}^{-1}(A(x,\xi)-A(m,n))\widehat{u}$ is really a reminder term of $A(x,D)u$ and one needs to find suitable function spaces to control the nonlinearity.   In the following we introduce some function spaces connected to the frequency-physical uniform decomposition operators.
For convenience, we write
\begin{align}\label{space1}
  \ell^q_s \ell^p(L^r):=\left\{\{g_{m,n}\}_{m,n\in\mathbb{Z}^d}:g_{m,n}\in\mathscr{S}'
  (\mathbb{R}^d), \left\| \| \{\langle n\rangle^s\|g_{m,n}\|_{r}\}   \|_{\ell^p_{m\in \mathbb{Z}^d}} \right\|_{\ell^q_{n\in \mathbb{Z}^d}}<\infty\right\}.
\end{align}
 In order to apply the frequency-physical-uniform decomposition operators, combining $\{\Box_{m,n}\}_{m,n\in \mathbb{Z}^d}$ with $\ell^q_s \ell^p(L^r) $, we can introduce the following norm:
\begin{align}\label{space0}
  \|f\|_{X_{r,p,q}^s}= \left\| \big\| \| \{\langle n\rangle^s\Box_{m,n} f\} \|_{r}  \big\|_{\ell^p_{m\in \mathbb{Z}^d}} \right\|_{\ell^q_{n\in \mathbb{Z}^d}}.
\end{align}
Combining \eqref{space0} and \eqref{space1}, we see that $\|f\|_{X_{r,p,q}^s}=\|\{\Box_{m,n}f\}\|_{\ell^q_s\ell^p(L^r)}.$  $X^s_{r,p,q}$ is of importance for us solving  \eqref{localization integral}. In fact, we can show that
\begin{prop}
Let $s\in \mathbb{R}$, $1\leq p,q,r\leq \infty$. Then
 $\|\cdot\|_{X^s_{r, p,q}}$  is an equivalent norm on $M^s_{p,q}$.
\end{prop}

In order to handle the function $u(x,t)$ of $(x,t) \in \mathbb{R}^d\times [0,T]$, we also need the following function spaces  $\mathscr{L}^\gamma(0,T;X^s_{r,p,q})$ for which the norm is defined as
\begin{align}\label{defn}
  \|u\|_{\mathscr{L}^\gamma(0,T;X^s_{r,p,q})}:=\Bigg(\sum_{n\in\mathbb{Z}^d}\langle n\rangle^{sq}\bigg(\sum_{m\in\mathbb{Z}^d}\|\Box_{m,n}u\|_{L^\gamma(0,T;L^r)}^p\bigg)^{q/p}\Bigg)^{1/q}<\infty.
\end{align}
Note that for $1\leq p\vee q\leq \gamma\leq \infty$, by Minkowski's inequality, we have $\mathscr{L}^\gamma(0,T;X^s_{r,p,q})\subset L^\gamma(0,T;X^s_{r,p,q})$. On the contrary, for $1\leq \gamma\leq p\wedge q\leq \infty$, we have $L^\gamma(0,T;X^s_{r,p,q})\subset \mathscr{L}^\gamma(0,T;X^s_{r,p,q})$.

In this paper, we will show the well posedness of \eqref{PDE} for initial data in   $M_{p,q}^s$ and $H^s$.  In order to indicate our method in a comparatively easier way, we first consider a simple case $A(x,D)= a(x)+b(D)$.  We assume that the following conditions are satisfied:
\begin{itemize}

\item[\rm (H1)]  $a(x)\in {C}^{N}(\mathbb{R}^d)$, $b(\xi)\in {C}^{N}(\mathbb{R}^d\setminus\{0\})$ for some $N\gg d$;

  \item[\rm (H2)] There exists $\sigma_1>0$ such that $|\partial_x^\alpha a(x)|\lesssim \langle x\rangle^{\sigma_1-|\alpha|}$ for $x\in \mathbb{R}^d$, \ $1\leq |\alpha| \leq N$;

  \item[\rm (H3)] There exists $\sigma_2>0$ such that $|\partial_\xi^\beta b(\xi)|\lesssim |\xi|^{\sigma_2-|\beta|}$, if $\xi\neq 0$ and $0\leq |\beta|\leq N$;

  \item[\rm (H4)] $\mathfrak{Re} a(x) \gtrsim |x|^{\sigma_1}$ and $\mathfrak{Re} b(\xi) \gtrsim |\xi|^{\sigma_2}$ for $x,\xi \in \mathbb{R}^d$;

\end{itemize}

First, we consider the general nonlinearity which contains derivative terms.

 \begin{thm} \label{thma}
 Let $1\leq p,q,r\leq \infty$.  Assume that the symbol $A(x,\xi)$ satisfies (H1)--(H4). Suppose that the nonlinearity takes the form as in \eqref{nonlinearity}, with $\kappa \leq\sigma_2$. We have the following results:

\begin{itemize}
\item[\rm (i)] Let $0<\kappa<\sigma_2$ and
  $$
    \gamma(K):= \max\left(K, \  \frac{\sigma_2(K-1)}{\sigma_2-\kappa}  \right), \ \ s(\rho):= \kappa+ \frac{d}{q'} - \frac{\sigma_2}{\rho}.
  $$
We assume that
$$
\gamma >
\left\{
\begin{array}{ll}
\gamma(K), &  if \ s(\gamma(K))\geq 0,\\
\sigma_1   \vee \gamma(K),  & if \   s(\gamma(K)) <0.
\end{array}
\right.
$$

Suppose that $s_0>\kappa+d/q'-\sigma_2/\gamma$,  $u_0\in M^{s_0}_{p,q}$.   Then there exists $T >0$ such that \eqref{localization integral} has a unique solution $u\in C([0,T]; M^{s_0}_{p,q})\cap \mathscr{L}^\gamma  (0,T; X^{s_0+\sigma_2/\gamma}_{r,p,q})$.

\item[\rm (ii)] Let $\kappa=\sigma_2$, $s_0>\sigma_2+d/q'$ and $u_0\in M^{s_0}_{p,q}$.  Suppose that $\|u_0\|_{M^{s_0}_{p,q}}$ is sufficiently small. Then  there exists $T>0$ such that \eqref{localization integral} has a unique solution $u\in C([0,T]; M^{s_0}_{p,q})\cap \mathscr{L}^\infty (0,T; X^{s_0}_{r,p,q})$.

\item[\rm (iii)] Let $0\leq \kappa\leq \sigma_2$.  Assume that $a(0)>0$ and $\mathfrak{ Re}\ (F(u),u)\leq 0$\footnote{We denote by $(\cdot, \ \cdot)$ the complex inner product.}.  Let $s_0> \kappa+d/2 $, $u_0\in H^{s_0}$ is sufficiently small.  Then \eqref{localization integral} has a unique global solution $u$ in $C(0,\infty; H^{s_0})\cap {L}^2(0,\infty;H^{s_0+\sigma_2/2})$.
\end{itemize}

\end{thm}

Next, if the nonlinearity contains no derivative terms, i.e., $\kappa=0$, we can assume that the initial data  have lower regularity. In particular, we can consider a class of data in modulation spaces with negative regularity index.

 \begin{thm} \label{thma0}
 Let $1\leq p,q, r\leq \infty$.  Assume that the symbol $A(x,\xi)$ satisfies (H1)--(H4). Suppose that the nonlinearity takes the form
$$
F(u)= \sum^I_{i=1} \lambda_i |u|^{k_i-1}u, \ \ k_i\in 2\mathbb{N}+1, \ \ K=\max_{1\leq i \leq I} k_i.
$$
We have the following results:
\begin{itemize}
\item[\rm (i)] Let $ \sigma_2 \leq d/q$, $s_0 \geq  0, \ s_0> d/q'-\sigma_2/(K-1)$, $K  \ll \gamma <\infty$ and $u_0\in M^{s_0}_{p,q}$. Then there exists $T >0$ such that \eqref{localization integral} has a unique solution $u\in C([0,T]; M^{s_0}_{p,q})\cap \mathscr{L}^\gamma  (0,T ; X^{s_0+\sigma_2/\gamma}_{r,p,q})$.

\item[\rm (ii)] Let $\sigma_2 >d/q$, $s\geq 0$ and $s >  d/q' - d/q(K-1)$. Denote
\begin{align} \label{gammaqKa}
    \gamma(q,K):= \max\left(K,  \ \frac{K-1}{1-d/q\sigma_2}  \right).
\end{align}
Assume that $s-\sigma_2/\gamma \geq 0$ for some $\gamma > \gamma(q,K)$, or  $s-\sigma_2/\gamma < 0$ for some $\gamma > \gamma(q,K) \vee \sigma_1$.
If $u_0\in M^{s-\sigma_2/\gamma}_{p,q}$, then there exists $T >0$ such that \eqref{localization integral} has a unique solution $u\in C([0,T]; M^{s-\sigma_2/\gamma}_{p,q})\cap \mathscr{L}^\gamma  (0,T ; X^{s }_{r, p,q})$.

\item[\rm (iii)] Let $d< 2\sigma_2$, $\mathfrak{Re}\lambda_i\leq 0$ for $1\leq i\leq I$, $u_0\in L^2$, and $K<1+ 2\sigma_2/d$. There exists $\gamma > \gamma (2,K)$ such that \eqref{localization integral} has a unique global solution $u$ in $C(0,\infty; L^2)\cap \mathscr{L}^\gamma_{\rm loc} (0,\infty; X^{\sigma_2/\gamma}_{2, 2,2})$.

\end{itemize}

\end{thm}

\begin{rem}\label{rem}
\rm Theorems \ref{thma} and \ref{thma0}  need several remarks.
\begin{itemize}

\item[(i)] $A(x,D)= \langle x\rangle^{\sigma_1}+ (-\Delta)^{\sigma_2/2}$ is included in Theorems \ref{thma} and \ref{thma0}.

\item[(ii)] Theorems \ref{thma} and \ref{thma0}  contain  a class of data in modulation spaces with negative regularity index. For example, we consider the case $q=1$ in (ii) of Theorem \ref{thma0}: If $\sigma_2>d$, $\gamma>\gamma(1,K) \vee \sigma_1$,  $u_0 \in M^{-\sigma_2/\gamma}_{p,1}$, then \eqref{localization integral} with the nonlinearity in Theorem \ref{thma0} is local well-posed in $C(0,T; M^{-\sigma_2/\gamma}_{p,1}) \cap \mathscr{L}^\gamma (0,T; X^0_{r,p,1})$ for some $T>0$.

\item[(iii)] Theorem \ref{thma} also contains the following semi-linear parabolic equation
\begin{align} \label{2ndPDE}
u_t + \langle x\rangle^{\sigma_1} u -\Delta u = -|\Delta u|^{K-1}  u, \ \ u(0,x)=u_0(x).
\end{align}
By the result of (iii) of Theorem \ref{thma}, we see that it is globally well-posed in the space $C(0,\infty; H^{2+d/2+\varepsilon}) \cap {L}^2(0,\infty;H^{3+d/2+\varepsilon})$  if $\sigma_1 >0$, $K\in 2\mathbb{N}+1$ and $u_0 \in H^{2+d/2+\varepsilon}$ is sufficiently small.

\item[(iv)] Recall that for the semi-group $e^{t\Delta}$,
$$
\|\Delta e^{t\Delta} u_0\|_p \lesssim t^{-1} \|u_0\|_p,
$$
one cannot obtain the global well-posedness of \eqref{2ndPDE} by only use the above decaying estimate, since $t^{-1}$ contains singularity at $t=0$.

\end{itemize}
\end{rem}
Finally, we consider the general case  $A\in S^M$. We have

\begin{thm}\label{thmc} Let $1\leq p,q,r\leq \infty$. Assume that the symbol $A\in S^M  $ is $M$-order dissipative and there exists $\varepsilon>0$ such that
$$
|\partial_x^\alpha A(x,\xi)|\leq A_\alpha(1+|\xi|)^{M-\varepsilon} \ for \   \ 1\leq |\alpha|\leq N.
$$
Suppose that the nonlinearity takes the form as in \eqref{nonlinearity}  with $\kappa \leq M$.  Let
  $$
  \gamma>
  \max\left(K, \frac{M(K-1)}{M-\kappa}  \right)\ \  { for} \ \kappa <M.
  $$
 We have the following results:
\begin{itemize}
\item[\rm (i)] Let $\kappa<M$, $s_0>\kappa+d/q'-M/\gamma$, $u_0\in M^{s_0}_{p,q}$.  Then there exists a positive time $T>0$ such that \eqref{localization integral} has a unique solution $u$ in $C([0,T]; M^{s_0}_{p,q})\cap \mathscr{L}^\gamma (0,T; X^{s_0+M/\gamma}_{r,p,q})$.

\item[\rm (ii)] Let $\kappa=M$, $s_0>M+d/q'$, $u_0\in M^{s_0}_{p,q}$ and $\|u_0\|_{M^{s_0}_{p,q}}$ is sufficiently small.  Then there exists a positive time $T$ such that \eqref{localization integral} has a unique solution  $u \in C(0,T; M^{s_0}_{p,q})\cap \mathscr{L}^\infty (0,T; X^{s_0}_{r,p,q})$.

\item[\rm (iii)] Assume that the symbol $A(0,0)>0$,  $\mathfrak{ Re}\ (F(u),u)\leq 0$,  and the pseudo-differential operator $A(x,D)$ satisfies
 $$\mathfrak{Re}\ (A(x,D)\varphi,\varphi)\gtrsim \|\varphi\|^2_{H^{M/2}}, \ \ \forall\ \varphi\in \mathscr{S}.$$
 Let $s_0>\kappa+d/2$, $u_0\in H^{s_0}$ and $\|u_0\|_{H^{s_0}}$ be sufficiently small, then  \eqref{localization integral} has a unique global solution $u$ in $C(0,\infty; H^{s_0})\cap {L}^2(0,\infty;H^{s_0+M/2})$.
\end{itemize}
\end{thm}

Some notations in this paper are as follows. For any multi-index $\alpha=(\alpha_1,\alpha_2,\cdots,\alpha_d)$, we denote $|\alpha|=\alpha_1+\alpha_2+\cdots+\alpha_d$,   $D_j^{\alpha_j}=({\rm - i})^{|\alpha_j|}\partial_j^{\alpha_j}$ for $j=1,2,\cdots,d$, and $D^\alpha= D_1^{\alpha_1}D_2^{\alpha_2}\cdots D_d^{\alpha_d}$. For any $p\in[1,\infty]$, $p'$ denotes the conjugate number of $p$, i.e., $1/p+1/p'=1$. We will use the Lebesgue spaces $L^p:=L^p(\mathbb{R}^d)$ with the norm $\|\cdot\|_p:=\|\cdot\|_{L^p(\mathbb{R}^d)}$. We denote by $\mathscr{S}:=\mathscr{S}(\mathbb{R}^d)$  the Schwartz space and by $\mathscr{S}':=\mathscr{S}'(\mathbb{R}^d)$ its dual space. we denote by $a\vee b= \max (a,b)$, $a\wedge b = \min(a,b)$, $a+= a+\varepsilon$ for some $0<\varepsilon\ll 1$. $\chi_E$ denotes the characteristic function on $E$.

The paper is organized as follows. In Section 2 we show $X^s_{r,p,q}=M^s_{p,q}$ with equivalent norms. We will obtain some algebraic properties of $\mathscr{L}^\gamma (0,T; X^{s}_{r,p,q})$ in Section 3. In Section 4 we consider the estimates of the linear terms in the integral equation \eqref{localization} in the spaces $\mathscr{L}^\gamma (0,T; X^{s}_{r,p,q})$.  In Section 5 we will make the estimates of the nonlinear term in the integral equation \eqref{localization} in the spaces $\mathscr{L}^\gamma (0,T; X^{s}_{r,p,q})$. The local and global well-posedness of \eqref{PDE} will be shown in Sections 6--8. In the Appendix, we show that the solution of \eqref{localization integral} solves \eqref{PDE} in some weak function spaces.

\section{Equivalent Norms}

First, let us recall the original definition of modulation spaces given by Feichtinger \cite{Fe83}.  The short-time
Fourier transform of a function $f$ with respect to $\varphi \in
\mathscr{S}$ is
\begin{align*}
V_{\varphi} f(x, \omega) = \int_{\mathbb{R}^d} e^{-{\rm i} t \cdot \omega}
\overline{\varphi (t-x)} f(t)dt,
\end{align*}
where $\varphi$ is said to be a window function. We write
\begin{align*}
\|f\|^{\circ}_{M^s_{p,q}} = \left(\int_{\mathbb{R}^n}\left(
\int_{\mathbb{R}^d} |V_{\varphi} f(x, \omega)|^p dx\right)^{q/p}
\langle\omega\rangle ^{sq} d\omega \right)^{1/q} 
\end{align*}
with the usual modifications if $p$ or $q$ is infinite.  Modulation spaces $M^s_{p,q}$ are defined as the
spaces of all distributions $f\in \mathscr{S}' $ for
which $\|f\|^{\circ}_{M^s_{p,q}}$ is finite (see Feichtinger
\cite{Fe83}).   A time--frequency localized definition  of $M^0_{1,1}$ was given by D\"orfler, Feichtinger and Gr\"ochenig \cite{DoFeGr06} via the localized operators
$$H_\lambda f: = V^*_{\varphi} \sigma_\lambda V_{\varphi} f, \ \ \lambda\in \mathbb{Z}^{2d}   $$
and they showed that $ \|\|H_{\lambda} \cdot\|_{L^2(\mathbb{R}^{2d})}  \|_{\ell_{1,1}}$ is  an equivalent norm on $M^0_{1,1}$. By using a different technique, D\"orfler and Gr\"ochenig \cite{DoGr11} generalized this equivalence to all modulation spaces $M^s_{p,q}$, i.e., $ \|\|H_{\lambda}  \cdot\|_{L^2(\mathbb{R}^{2d})}  \|_{\ell^s_{p,q}}$ and $\|\cdot\|_{M^s_{p,q}}$ are equivalent norms for all $1\leq p,q \leq \infty$.

Using the operators $\Box_{m,n}$,  in this section we show that $M^s_{p,q}= X^s_{r,p,q}$  with  equivalent norms.  Noticing that $\Box_n f$  can be regarded as a frequency-localized version of $V_{\varphi} f$, one sees that $\Box_{m,n} f$ is a frequency-physical localized version of $V_{\varphi} f$, which  is more closely related  to $\sigma_m \sigma_n V_{\varphi} f$ than $H_{(m,n)} f$.  We start with the following lemma, which is also useful in the next few sections.

\begin{lem}\label{lemma1}  Let $1\leq p,q\leq\infty$. Assume that $\varphi,\psi_k $ are Schwartz functions satisfying
${\rm supp} \psi_k\subset\{x:\ |x-k|_\infty\leq1\}$ and
\begin{gather}
    \|\psi_k\|_{p'}+\|\psi_k\|_{q} \leq C_k, \ \ k\in \mathbb{Z}^d,  \nonumber\\
    C_{N,\varphi}:=\sup_{|\alpha|\leq N}\|\partial^\alpha\varphi\|_1 , \ \ N\in\mathbb{N}. \nonumber
\end{gather}
Denote for $m,k\in \mathbb{Z}^d$,
\begin{align*}
T_{m,k}f:=\psi_m(x) \mathscr{F}^{-1} \varphi(\xi) \mathscr{F}\psi_k(x)f.
\end{align*}
Then we have
\begin{align*}
\|T_{m,k}f\|_q\lesssim C_m C_k C_{N,\varphi}\langle m-k\rangle^{-N}\|f\|_p.
\end{align*}
\end{lem}
{\bf Proof.} Let us rewrite
\begin{align}
T_{m,k}f=\psi_m(x) \int_{\mathbb{R}^d}\psi_k(y)f(y)\int_{\mathbb{R}^d} e^{{\rm i}(x-y)\cdot \xi} \varphi(\xi) d\xi dy.\label{lemma1e}
\end{align}
Obviously,
\begin{align}
\bigg|\int_{\mathbb{R}^d} e^{{\rm i}(x-y)\cdot \xi} \varphi(\xi) d\xi \bigg|\leq\|\varphi\|_1\leq C_{N,\varphi}.\label{lemma1f}
\end{align}
Integrating by parts, we have
\begin{align*}
\bigg|\int_{\mathbb{R}^d} e^{{\rm i}(x-y)\cdot \xi} \varphi(\xi) d\xi \bigg|=\bigg|\frac{1}{(x-y)^\alpha}\int_{\mathbb{R}^d} e^{{\rm i}(x-y)\cdot \xi} \partial^\alpha\varphi(\xi) d\xi\bigg|.
\end{align*}
Hence, we have
\begin{align}
\bigg|\int_{\mathbb{R}^d} e^{{\rm i}(x-y)\cdot \xi} \varphi(\xi) d\xi \bigg|\leq C_{N,\varphi}\min(1,|x-y|^{-N})\leq (2d)^N C_{N,\varphi}\langle x-y\rangle^{-N}.\label{lemma1h}
\end{align}
Inserting the estimate of \eqref{lemma1h} into \eqref{lemma1e}, we have
\begin{align*}
|T_{m,k}f(x)|\leq\bigg|\psi_m(x) \int_{\mathbb{R}^d} C^N C_{N,\varphi}\langle x-y\rangle^{-N}|\psi_k(y)f(y)| dy\bigg|.
\end{align*}
By translation and H\"older's inequality,
\begin{align*}
\|T_{m,k}f\|_q&= \left\|\psi_m(\cdot+m)\int_{\mathbb{R}^d} C^N C_{N,\varphi}\langle m-k+ x-y\rangle^{-N}|\psi_k(y+k)f(y+k)| dy \right\|_q\nonumber\\
&\leq C^N C_{N,\varphi}\langle m-k\rangle^{-N}\|\psi_m\|_q \|\psi_k\|_{p'}\|f\|_p\nonumber\\
&\leq   C^N C_m C_k C_{N,\varphi}\langle m-k\rangle^{-N}\|f\|_p.
\end{align*}
The result follows. $\hfill\Box$\\

{\bf Remark.} The result in Lemma \ref{lemma1} is also right for non-integral order decay. Indeed, combining \eqref{lemma1f} and \eqref{lemma1h}, for $\theta\in[0,1]$, we have
\begin{align}\label{remark1}
\|T_{m,k}f\|_q\lesssim C_m C_k C_{N,\varphi}\langle m-k\rangle^{-\theta N}\|f\|_p.
\end{align}

\begin{prop}\label{embedding} {\rm (Equivalent Norm)}
Let $s \in \mathbb{R}$ and $1\leq   p,q,  r \leq\infty$. Then $X_{r, p,q }^{s }  = M^s_{p,q}$ with equivalent norm.  In particular, we have $X^s_{r,2,2} =H^s$ with equivalent norm.
\end{prop}
{\bf Proof.} First, we show that
$$
X_{r_1,p,q}^{s}\subset X_{r_2,p,q}^{s}, \ \  \forall \ 1\leq r_1,  r_2\leq \infty.
$$
Let us observe that
$$
\Box_{m,n}f=\sigma_m(x) \mathscr{F}^{-1} \sigma_n(\xi) \mathscr{F}f=\sum_{|l_1|_\infty \vee |l_2|_\infty\leq1}\sum_{k\in \mathbb{Z}^d}\sigma_m(x) \mathscr{F}^{-1} \sigma_{n+l_1}(\xi) \mathscr{F}\sigma_{k+l_2}(x)\Box_{k,n}f.
$$
From Lemma \ref{lemma1} it follows that
\begin{align}\label{embedding0}
\|\Box_{m,n}f\|_{r_2}\lesssim\sum_{k\in \mathbb{Z}^d}\langle m-k\rangle^{-2d}\|\Box_{k,n}f\|_{r_1}.
\end{align}
By definition,
\begin{align}\label{embeddinga}
\|f\|_{X_{r_2,p,q}^{s}}=\Bigg(\sum_{ n\in\mathbb{Z}^d} \langle n\rangle^{sq}  \bigg( \sum_{ m \in\mathbb{Z}^d} \|\Box_{m,n}f\|_{r_2}^{p}  \bigg)^{q/p}         \Bigg)^{1/q}.
\end{align}
Inserting \eqref{embedding0} into \eqref{embeddinga}, and using Young's inequality, we have
\begin{align*}
\|f\|_{X_{r_2,p, q}^{s} }& \lesssim
\Bigg(\sum_{n\in\mathbb{Z}^d}\langle n\rangle^{sq}   \bigg( \sum_{m\in\mathbb{Z}^d}\bigg( \sum_{k\in \mathbb{Z}^d}\langle m-k\rangle^{-2d}\|\Box_{k,n}f\|_{r_1} \bigg)^p  \bigg)^{q/p} \Bigg)^{1/q}\nonumber\\
& \lesssim \Bigg(\sum_{n\in\mathbb{Z}^d} \langle n\rangle^{s q } \bigg(\sum_{k\in \mathbb{Z}^d} \|\Box_{k,n}f\|^{p}_{r_1} \bigg)^{q/p} \Bigg)^{1/q}=\|f\|_{X_{r_1,p, q}^{s}}.
\end{align*}
 The embedding $X_{r_1,p,q}^{s}\subset X_{r_2,p,q}^{s}$ is obtained.  Exchanging the roles of $r_1$ and $r_2$, we have $X_{r_2,p,q}^{s} \subset X_{r_1,p,q}^{s}$. Then we have  $X_{r_1,p,q}^{s}= X_{r_2,p,q}^{s}$.
It follows that  $X_{r, p,q}^{s} =X^s_{p,p,q}$. Let us observe
 \begin{align*}
\|f\|_{X_{p,p,q}^{s}}=\Bigg(\sum_{n\in\mathbb{Z}^d}\langle n\rangle^{s q} \bigg(\int \Big(\sum_{m\in \mathbb{Z}^d} |\sigma_m(x)|^p \Big) |\Box_{n}f (x)|^{p} dx\bigg)^{q/p} \Bigg)^{1/q}.
\end{align*}
Noticing that $\sum_{m\in \mathbb{Z}^d} |\sigma_m(x)|^q =\sum_{|l|_\infty \leq 1} |\sigma_{m+l}(x)|^q \sim 1$ for any $x\in \mathbb{R}^d$, we immediately have $X_{p,p,q}^{s}=M^s_{p,q}$.   $\hfill\Box$\\

{\bf Remark.} One may ask the question why $X_{r, p,q}^{s}$ is independent of $r\in [1,\infty]$. It is well-known that for $r_1\leq r_2$,  $\|f\|_{r_2} \lesssim \|f\|_{r_1}$ if ${\rm supp} \widehat{f}$ is compact,  $\|f\|_{r_1} \lesssim \|f\|_{r_2}$ if  ${\rm supp} {f}$ is compact. Hence,  $\|f\|_{r_1}  \sim  \|f\|_{r_2}$ if $ f$ and $ \widehat{f}$ could be both supported in compact sets. Unfortunately, it is impossible to expect that ${\rm supp} f$ and ${\rm supp} \widehat{f}$ are both compact.  By a formal observation, we see that the support set of $\Box_{m,n}f$ is compact in physical space and $\mathscr{F}^{-1}\sigma_n\mathscr{F}f$ has a compact support set in frequency space. Even though we cannot realize
\begin{align}\label{embeddingrmk}
\|\Box_{m,n}f\|_{r_2}\lesssim \|\Box_{m,n}f\|_{r_1}, \qquad  r_2\geq r_1,
\end{align}
we can show  \eqref{embedding0} which is an approximative version of \eqref{embeddingrmk}.  $\|\Box_{m,n}f\|_{r_2}$  is mainly bounded by  $\|\Box_{k,n}f\|_{r_1}$ for $|k-m|\lesssim 1$, and the contributions of the other $\|\Box_{k,n}f\|_{r_1}$ are very small due to the fast decay of $\langle m-k\rangle^{-N}$.

\section{Multi-linear Estimates}

\begin{prop}\label{bilinear} {\rm (Algebraic Estimates)}
Let   $1\leq p,q,r\leq\infty$. We have
\begin{itemize}
\item[\rm (i)] Let $s\geq 0$. Then $X^s_{r,p,1}$ is a Banach algebra:
\begin{align}\label{bilineara}\|u_1u_2\|_{X^s_{r,p,1}}\lesssim \|u_1\|_{X^s_{r,p,1}}\|u_2\|_{X^s_{r,p,1}}.
\end{align}
\item[\rm (ii)] Let $q>1$, $s> d/q'$. Then $X^s_{r,p,q}$ is a Banach algebra:
\begin{align}\label{bilinearb}
\|u_1u_2\|_{X^s_{r,p,q}}\lesssim \|u_1\|_{X^s_{r,p,q}}\|u_2\|_{X^s_{r,p,q}}.
\end{align}
\end{itemize}
\end{prop}
{\bf Proof.} In view of $X_{r,p,q}^{s }=M^{s}_{p,q}$  and   Proposition 6.9 and Remark 6.4 in \cite{Fe83}, we have the results. $\hfill\Box$

Now we consider the algebraic structure of  $\mathscr{L}^\gamma(0,T;X^s_{r,p,q})$.

\begin{prop}\label{nonlinear}{\rm (Multi-linear Estimates)}
Let $1\leq p,q,r,\gamma,p_i,q_i,r_i,\gamma_i\leq\infty$ satisfy
\begin{align*}
\frac{1}{\gamma}=\frac{1}{\gamma_1}+\frac{1}{\gamma_2}+\cdots+\frac{1}{\gamma_N},\ \ \ \  &\frac{1}{r}=\frac{1}{r_1}+\frac{1}{r_2}+\cdots+\frac{1}{r_N},\\
\frac{1}{p}=\frac{1}{p_1}+\frac{1}{p_2}+\cdots+\frac{1}{p_N},\ \ \ \ & \frac{1}{q}=\frac{1}{q_1}+\frac{1}{q_2}+\cdots+\frac{1}{q_N}-(N-1),
\end{align*}
 and $s\geq 0$. Then we have
\begin{align}\label{nonlinearA}
\|f_1\cdots f_N\|_{\mathscr{L}^\gamma(0,T;X^s_{r,p,q})}\lesssim &  \|f_1\|_{\mathscr{L}^{\gamma_1}(0,T;X^{s}_{r_1,p_1,q_1})} \prod^N_{i=2}\|f_i\|_{\mathscr{L}^{\gamma_i}(0,T;X^{0}_{r_i,p_i,q_i})} +... \nonumber\\
&  + \prod^{N-1}_{i=1}\|f_i\|_{\mathscr{L}^{\gamma_{i+1}}(0,T;X^{0}_{r_{i+1},p_{i+1},q_{i+1}})} \|f_N\|_{\mathscr{L}^{\gamma_1}(0,T;X^{s}_{r_1,p_1,q_1})}.
\end{align}
In particular, if $q=1$, we have
\begin{align}
\|f_1 f_2\cdots f_N\|_{\mathscr{L}^\gamma(0,T;X^s_{r,p,1})}\lesssim & \|f_1\|_{\mathscr{L}^{\gamma_1}(0,T;X^{s}_{r,p,1})} \prod^N_{i= 2}\|f_i\|_{\mathscr{L}^{\gamma_i}(0,T;X^{0}_{r,p,1})}
+... \nonumber\\
& +  \prod^{N-1}_{i= 1}\|f_i\|_{\mathscr{L}^{\gamma_{i+1}}(0,T;X^{0}_{r,p,1})}
\|f_N\|_{\mathscr{L}^{\gamma_1}(0,T;X^{s}_{r,p,1})}
.\label{nonlinearB}
\end{align}
If $q>1$, we have
\begin{align}
\|f_1 f_2\cdots f_N\|_{\mathscr{L}^\gamma(0,T;X^s_{r,p,q})}\lesssim &  \|f_1\|_{\mathscr{L}^{\gamma_1}(0,T;X^{s}_{r,p,q})} \prod^N_{i=2 }\|f_i\|_{\mathscr{L}^{\gamma_i}(0,T;X^{d/q'+}_{r,p,q})} +...\nonumber\\
&  + \|f_N\|_{\mathscr{L}^{\gamma_1}(0,T;X^{s}_{r,p,q})} \prod^{N-1}_{i=1 }\|f_{i}\|_{\mathscr{L}^{\gamma_{i+1}}(0,T;X^{d/q'+}_{r,p,q})}
. \label{nonlinearC}
\end{align}
\end{prop}
{\bf Proof.} By definition
\begin{align*}
\|f_1 f_2\cdots f_N\|_{\mathscr{L}^\gamma(0,T;X^s_{r,p,q})}=\bigg\|\langle n\rangle^s\Big\|\|\Box_{m,n}(f_1 f_2\cdots f_N)\|_{L^\gamma(0,T;L^r)}\Big\|_{\ell^p_m}\bigg\|_{\ell^q_n}.
\end{align*}
It is easy to see that
\begin{align*}
\|\Box_{m,n}(f_1 f_2\cdots f_N)\|_{L^\gamma(0,T;L^r)}\leq\sum_{n_1,n_2,\cdots,n_N\in\mathbb{Z}^d}\big\|\Box_{m,n}\big(\Box_{n_1}f_1\Box_{n_2}f_2\cdots\Box_{n_N}f_N\big)\big\|_{L^\gamma(0,T;L^r)}.
\end{align*}
From $\Box_n\big(\Box_{n_1}f_1\ \Box_{n_2}f_2\cdots\Box_{n_N}f_N\big)=0$ if $|n-n_1-\cdots-n_N|\gg 1$, it follows  that
\begin{align*}
&\|\Box_{m,n}(f_1 f_2\cdots f_N)\|_{L^\gamma(0,T;L^r)}\nonumber\\
&\leq\sum_{n_1,n_2,\cdots,n_N\in\mathbb{Z}^d}\big\|\Box_{m,n}\big(\Box_{n_1}f_1\Box_{n_2}f_2\cdots\Box_{n_N}f_N\big)\big\|_{L^\gamma(0,T;L^r)} \chi_{|n-n_1-\cdots-n_N|\lesssim 1}.
\end{align*}
Denote
\begin{align} \label{box}
\widetilde{\Box}_{m,n}:= \sum_{|h|_\infty\leq 1}\sigma_{m+h}(x) \mathscr{F}^{-1} \sigma_n(\xi) \mathscr{F},
\end{align}
and from Lemma \ref{lemma1}, we have
\begin{align*}
&\|\Box_{m,n}(f_1 f_2\cdots f_N)\|_{L^\gamma(0,T;L^r)}\nonumber\\
&\leq\sum_{n_1,n_2,\cdots,n_N\in\mathbb{Z}^d}\sum_{k\in\mathbb{Z}^d}\big\|\Box_{m,n}\big(\Box_{k,n_1}f_1
\widetilde{\Box}_{k,n_2}f_2\cdots\widetilde{\Box}_{k,n_N}f_N\big)\big\|_{L^\gamma(0,T;L^r)} \chi_{|n-n_1-\cdots-n_N|\lesssim 1}\nonumber\\
&\lesssim \sum_{k\in\mathbb{Z}^d}\langle m-k\rangle^{-2d} \sum_{n_1,n_2,\cdots,n_N\in\mathbb{Z}^d}
\big\|\Box_{k,n_1}f_1\widetilde{\Box}_{k,n_2}f_2\cdots\widetilde{\Box}_{k,n_N}f_N\big\|_{L^\gamma(0,T;L^r)} \chi_{|n-n_1-\cdots-n_N|\lesssim 1}.
\end{align*}
By H$\ddot{\rm o}$lder's inequality, we get that
\begin{align*}
&\|\Box_{m,n}(f_1 f_2\cdots f_N)\|_{L^\gamma(0,T;L^r)}\nonumber\\
&\lesssim \sum_{k\in\mathbb{Z}^d}\langle m-k\rangle^{-2d} \sum_{n_1,n_2,\cdots,n_N\in\mathbb{Z}^d}\|\Box_{k,n_1}f_1\|_{L^{\gamma_1}(0,T;L^{r_1})}
\prod^N_{i=2}\|\widetilde{\Box}_{k,n_i}f_i\|_{L^{\gamma_i}(0,T;L^{r_i})} \chi_{|n-n_1-\cdots-n_N|\lesssim 1}.
\end{align*}
From Minkowski's inequality, Young's inequality, and H$\ddot{\rm o}$lder's inequality,  we obtain
\begin{align}\label{nonlinear1}
&\Big\|\|\Box_{m,n}(f_1 f_2\cdots f_N)\|_{L^\gamma(0,T;L^r)}\Big\|_{\ell^p_m}\nonumber\\
&\lesssim \sum_{n_1,n_2,\cdots,n_N\in\mathbb{Z}^d}\bigg\|\|\Box_{m,n_1}f_1\|_{L^{\gamma_1}(0,T;L^{r_1})}
\prod^N_{i=2}\|\widetilde{\Box}_{m,n_i}f_i\|_{L^{\gamma_i}(0,T;L^{r_i})}\bigg\|_{\ell^p_m} \chi_{|n-n_1-\cdots-n_N|\lesssim 1}\nonumber\\
&\lesssim \sum_{n_1,n_2,\cdots,n_N\in\mathbb{Z}^d}\Big\|\|\Box_{m,n_1}f_1\|_{L^{\gamma_1}(0,T;L^{r_1})}\Big\|_{\ell^{p_1}_m}
\prod^N_{i=2}\Big\|\|\widetilde{\Box}_{m,n_i}f_i\|_{L^{\gamma_i}(0,T;L^{r_i})}\Big\|_{\ell^{p_i}_m} \chi_{|n-n_1-\cdots-n_N|\lesssim 1}.
\end{align}
Without loss of generality, we can assume that $\langle n_1\rangle=\max_{1\leq i\leq N} \langle n_i\rangle $. It follows from $|n-n_1-\cdots-n_N|\lesssim 1$ that $\langle n\rangle\lesssim \langle n_1\rangle$. Taking the $\ell^q_s$ norm in both sides of \eqref{nonlinear1} and using Young's inequality, we have from $\langle n\rangle^s\lesssim \langle n_1\rangle^s$,
\begin{align}\label{nonlinear2}
\|f_1 f_2\cdots f_N\|_{\mathscr{L}^\gamma(0,T;X^s_{r,p,q})}\lesssim \|f_1\|_{\mathscr{L}^{\gamma_1}(0,T;X^{s}_{r_1,p_1,q_1})}
\prod^N_{i=2}\|f_i\|_{\mathscr{L}^{\gamma_i}(0,T;X^{0}_{r_i,p_i,q_i})},
\end{align}
which implies \eqref{nonlinearA}. Observing  $X_{r_i,p,q }^{s}=X_{r,p,q}^{s}= M^s_{p,q}$ in Proposition \ref{embedding} and $\ell^p\subset \ell^{p_i}$ for $i\in\{1,2,\cdots,N\}$ and taking $q_1=q$, $q_i=1$ for $2\leq i\leq N$, we have from \eqref{nonlinear2} that
\begin{align}\label{nonlinear3}
\|f_1 f_2\cdots f_N\|_{\mathscr{L}^\gamma(0,T;X^s_{r,p,q})}\lesssim \|f_1\|_{\mathscr{L}^{\gamma_1}(0,T;X^{s}_{r,p,q})}
\prod^N_{i=2}\|f_i\|_{\mathscr{L}^{\gamma_i}(0,T;X^{0}_{r,p,1})}.
\end{align}
If $q=1$, \eqref{nonlinear3} implies \eqref{nonlinearB}. If $q>1$, from H$\ddot{\rm o}$lder's inequality,
\begin{align}\label{nonlinear4}
\|f_i\|_{\mathscr{L}^{\gamma_i}(0,T;X^{0}_{r,p,1})}&=\sum_{n\in\mathbb{Z}^d}\langle n\rangle^{-(d/q'+)}\langle n\rangle^{d/q'+}\Big\|\|\Box_{m,n}f_i\|_{L^{\gamma_i}(0,T;L^r)}\Big\|_{\ell^p_m}\nonumber\\
&\lesssim \bigg(\sum_{n\in\mathbb{Z}^d}\langle n\rangle^{-(d/q'+)q'}\bigg)^{1/q'}\|f_i\|_{\mathscr{L}^{\gamma_i}(0,T;X^{d/q'+}_{r,p,q})}\nonumber\\
&\lesssim \|f_i\|_{\mathscr{L}^{\gamma_i}(0,T;X^{d/q'+}_{r,p,q})}
\end{align}
Inserting \eqref{nonlinear4} into \eqref{nonlinear3}, we get the conclusion \eqref{nonlinearC}. $\hfill\Box$

\section{Linear Estimates for $A(x,D)=a(x)+b(D)$}

In this section, we will always assume that the symbol $A(x,\xi)=a(x) + b(\xi)$ satisfies the hypothesis (H1)--(H4). For convenience, we denote $s_\gamma:=s_0+\sigma_2/\gamma$, in particular, we have $s_\infty=s_0$, $s_2=s_0+\sigma_2/2$. We also write $\mathscr{L}^\gamma(0,T;\ell^q_{s} \ell^p(L^r)):=\ell^q_s\ell^p(L^\gamma(0,T; L^r))$.

\begin{lem}
\label{locallemma1} Let $1\leq p,q,r,\gamma\leq\infty$, $0<T<1$. Assume that $A(x,\xi)$ satisfies (H4).  Then we have
\begin{align*}
   \|\{e^{-t A(m,n)}\Box_{m,n} u_0\}\|_{\mathscr{L}^\gamma(0,T;\ell^q_{s_\gamma} \ell^p(L^r))}\lesssim \|u_0\|_{M^{s_0}_{p,q}}.
\end{align*}
\end{lem}
{\bf Proof.} Using condition (H4),
\begin{align}
\|e^{-t A(m,n)}\Box_{m,n} u_0\|_{L^\gamma(0,T;L^r)}&\lesssim\bigg(\int_0^T e^{-ct|n|^{\sigma_2}\gamma } dt\bigg)^{1/\gamma} \|\Box_{m,n} u_0\|_{L^r}\nonumber\\
&\lesssim \langle n\rangle^{-\sigma_2/\gamma}\|\Box_{m,n} u_0\|_{L^r}. \label{semigroupest1}
\end{align}
By the definition \eqref{space1} and \eqref{space0}, we get the desired conclusion. $\hfill\Box$

\begin{lem}\label{locallemma2}   Let  $1\leq p,q,r\leq\infty$, $1 <\gamma\leq\infty$, $s_\gamma\geq 0$.  In addition we assume that $\gamma>\sigma_1$ if $s_0<0$. Assume that $a(x)$ satisfies (H1),(H2) and (H4) for   $N= [s_\gamma]+d+2$. Then there exists $\delta>0$ such that
\begin{align}
& \|\mathscr{A}_{m,n}((a(x)-a(m))u)\|_{\mathscr{L}^\gamma(0,T;\ell^q_{s_\gamma} \ell^p(L^r))\cap \mathscr{L}^\infty(0,T;\ell^q_{s_0} \ell^p(L^r))}\lesssim T^\delta\|u\|_{\mathscr{L}^\gamma(0,T;X^{s_\gamma}_{r,p,q})}.\label{locallemma2a}\\
& \|\mathscr{A}_{m,n}((a(x)-a(m))u)\|_{{L}^\infty(0,T;\ell^q_{s_0} \ell^p(L^r))}\lesssim T^\delta\|u\|_{ {L}^\infty(0,T;X^{s_0}_{r,p,q})}, \ if \ s_0\geq 0.\label{locallemma2aaa}
\end{align}
\end{lem}
To prove this lemma,  we need the following lemma.
\begin{lem}\label{lemma4} Let $1\leq r\leq \infty$. Assume that $a(x)$ satisfies $|\partial_x^\alpha a(x)|\leq {{A_\alpha}} \langle x\rangle^{\sigma_1-|\alpha|}$ for any $1\leq |\alpha|\leq K+1$. Then for any $N>d+\sigma_1-1$, we have
\begin{align}\label{lemma4a}
&\|\Box_{m,n}((a(x)-a(m))u)\|_r\nonumber\\
&\lesssim \langle m\rangle^{(\sigma_1-1)\vee0}\sum_{n_1,n_2,l\in\mathbb{Z}^d}\langle m-l\rangle^{-N}\langle n_1\rangle^{-K}\|\widetilde{\Box}_{l,n_2}u\|_r \chi_{|n-n_1-n_2|_\infty \leq 3},
\end{align}
where $\widetilde{\Box}_{l,n} $ is as in \eqref{box}.
\end{lem}
{\bf Proof.}
By using $\sum_{k\in\mathbb{Z}^d} \sigma_{k}(x)=1$,  one can get
\begin{align}
\|\Box_{m,n}((a(x)-a(m))u)\|_r&\leq \sum_{k\in\mathbb{Z}^d}\|\Box_{m,n}(\sigma_k(x)(a(x)-a(m))u)\|_r.\nonumber
\end{align}
Noticing that $\Box_n\big(\Box_{n_1}u_1\ \Box_{n_2}u_2\big)=0$ if $|n-n_1-n_2|_\infty \geq 4$, we have
\begin{align*}
&\|\Box_{m,n}((a(x)-a(m))u)\|_r&\nonumber\\
&\leq\sum_{n_1,n_2,k\in\mathbb{Z}^d}\|\Box_{m,n}(\Box_{n_1}(\sigma_k(x)(a(x)-a(m))) \ \Box_{n_2}u)\|_r \chi_{|n-n_1-n_2| _\infty \leq 3}.
\end{align*}
Using the almost orthogonality $\sigma_l=\sum_{|h|_\infty\leq1}\sigma_l\sigma_{l+h}$,  from Lemma \ref{lemma1} it follows that for any $N\in\mathbb{N}$,
\begin{align}\label{lemma4b}
&\|\Box_{m,n}((a(x)-a(m))u)\|_r&\nonumber\\
&\leq\sum_{n_1,n_2,k,l\in\mathbb{Z}^d}\|\Box_{m,n}(\Box_{l,n_1}(\sigma_k(x)(a(x)-a(m))) \ \widetilde{\Box}_{l,n_2}u)\|_r \chi_{|n-n_1-n_2|_\infty \leq 3 }\nonumber\\
&\lesssim\sum_{n_1,n_2,k,l\in\mathbb{Z}^d}\langle m-l\rangle^{-N-1}\|\Box_{l,n_1} (\sigma_k(x)(a(x)-a(m))) \ \widetilde{\Box}_{l,n_2}u\|_r \chi_{|n-n_1-n_2|_\infty \leq 3}.
\end{align}
By H$\ddot{\rm o}$lder's inequality, \eqref{lemma4b} implies that
\begin{align}\label{lemma4e}
&\|\Box_{m,n}((a(x)-a(m))u)\|_r&\nonumber\\
&\lesssim\sum_{n_1,n_2,k,l\in\mathbb{Z}^d}\langle m-l\rangle^{-N-1}\|\Box_{l,n_1} (\sigma_k(x)(a(x)-a(m)))\|_\infty\|\widetilde{\Box}_{l,n_2}u\|_r \chi_{|n-n_1-n_2|_\infty \leq 3 }.
\end{align}
Using the multiplier estimates and then applying Lemma \ref{lemma1}, we immediately have for any $K,M\in\mathbb{N}$,
\begin{align}\label{lemma4f}
&\|\Box_{l,n_1} (\sigma_k(x)(a(x)-a(m)))\|_\infty\nonumber\\
&\lesssim \langle n_1\rangle^{-2K}\|\sigma_l(x)\mathscr{F}^{-1}\langle n_1\rangle^{2K}\langle \xi\rangle^{-2K}\sigma_{n_1}(\xi)\mathscr{F}(I-\Delta)^{K}(\sigma_k(x)(a(x)-a(m)))\|_\infty\nonumber\\
&\lesssim\langle n_1\rangle^{-2K}\langle l-k\rangle^{-M-1}\|(I-\Delta)^{K}(\sigma_k(x)(a(x)-a(m)))\|_\infty.
\end{align}
In view of the condition (UD), the mean value theorem and the assumption of $a(x)$, we get that
\begin{align}\label{lemma4g}
\|(I-\Delta)^{K} (\sigma_k(x)(a(x)-a(m)))\|_\infty&\lesssim\big(\langle k\rangle^{\sigma_1-1}+\langle m\rangle^{\sigma_1-1}\big)\langle m-k \rangle\nonumber\\
&\lesssim\big(\langle k\rangle^{\sigma_1-1}+\langle m\rangle^{\sigma_1-1}\big)\langle m-l\rangle \langle l-k\rangle.
\end{align}
Inserting \eqref{lemma4g} into \eqref{lemma4f} and noticing that $\langle n_1\rangle^{-2K} \lesssim \langle n_1\rangle^{-K}$, we have
\begin{align}\label{lemma4j}
\|\Box_{l,n_1} (\sigma_k(x)(a(x)-a(m)))\|_\infty\lesssim
\langle n_1\rangle^{-K}\langle l-k\rangle^{-M}\big(\langle k\rangle^{\sigma_1-1}+\langle m\rangle^{\sigma_1-1}\big)\langle m-l\rangle,
\end{align}
Therefore, inserting \eqref{lemma4j} into \eqref{lemma4e}, we have
\begin{align}\label{lemma4h}
&\|\Box_{m,n}((a(x)-a(m))u)\|_r&\nonumber\\
&\lesssim\sum_{n_1,n_2,k,l\in\mathbb{Z}^d}\langle m-l\rangle^{-N}\langle l-k\rangle^{-M}\big(\langle k\rangle^{\sigma_1-1}+\langle m\rangle^{\sigma_1-1}\big)\langle n_1\rangle^{-K}\|\widetilde{\Box}_{l,n_2}u\|_r \chi_{|n-n_1-n_2|_\infty \leq 3}.
\end{align}
If $0\leq\sigma_1\leq1$, then $\langle k\rangle^{\sigma_1-1}+\langle m\rangle^{\sigma_1-1}\lesssim1$. It follows that for $M>d$,
\begin{align*}
&\|\Box_{m,n}((a(x)-a(m))u)\|_r&\nonumber\\
&\lesssim\sum_{n_1,n_2,l\in\mathbb{Z}^d}\langle m-l\rangle^{-N}\bigg(\sum_{k\in\mathbb{Z}^d}\langle l-k\rangle^{-M}\bigg)\langle n_1\rangle^{-K}\|\widetilde{\Box}_{l,n_2}u\|_r \chi_{|n-n_1-n_2|_\infty \leq 3 }\nonumber\\
&\lesssim\sum_{n_1,n_2,l\in\mathbb{Z}^d}\langle m-l\rangle^{-N}\langle n_1\rangle^{-K}\|\widetilde{\Box}_{l,n_2}u\|_r \chi_{|n-n_1-n_2|_\infty \leq 3 }.
\end{align*}
If $\sigma_1>1$, we divide \eqref{lemma4h} into $\langle k\rangle\leq \langle m\rangle$ and $\langle k\rangle > \langle m\rangle$ two parts, which are denoted by ${\rm I}_{m,n}$ and ${\rm II}_{m,n}$ respectively. For $M>d$, we obtain that
\begin{align*}
{\rm I}_{m,n}
&\lesssim \langle m\rangle^{\sigma_1-1}\sum_{n_1,n_2,l\in\mathbb{Z}^d}\langle m-l\rangle^{-N}\bigg(\sum_{k\in\mathbb{Z}^d}\langle l-k\rangle^{-M}\bigg)\langle n_1\rangle^{-K}\|\widetilde{\Box}_{l,n_2}u\|_r \chi_{|n-n_1-n_2|_\infty \leq 3}\nonumber\\
&\lesssim\langle m\rangle^{\sigma_1-1}\sum_{n_1,n_2,l\in\mathbb{Z}^d}\langle m-l\rangle^{-N}\langle n_1\rangle^{-K}\|\widetilde{\Box}_{l,n_2}u\|_r \chi_{|n-n_1-n_2|_\infty \leq 3 }.
\end{align*}
Furthermore, we divide ${\rm II}_{m,n}$ into three parts, that is
\begin{align*}
{\rm II}_{m,n}&\lesssim \sum_{n_1,n_2,k,l\in\mathbb{Z}^d}\langle m-l\rangle^{-N}\langle l-k\rangle^{-M}\langle k\rangle^{\sigma_1-1}\langle n_1\rangle^{-K}\|\widetilde{\Box}_{l,n_2}u\|_r \chi_{|n-n_1-n_2|_\infty \leq 3 }\chi_{\{\langle k\rangle > \langle m\rangle\}}\nonumber\\
&\quad\quad\quad\quad\quad \times \Big(\chi_{\{\langle k\rangle \geq 10\langle l\rangle\}}+\chi_{\{\langle k\rangle \leq \frac{\langle l\rangle}{10}\}} +\chi_{\{\langle k\rangle \in(\frac{\langle l\rangle}{10}, 10\langle l\rangle)\}}\Big)\nonumber\\
&=:{\rm II}^{(1)}_{m,n}+{\rm II}^{(2)}_{m,n}+{\rm II}^{(3)}_{m,n}.
\end{align*}
For $\langle k\rangle \geq 10\langle l\rangle$, we have $\langle l-k\rangle^{-M}\lesssim \langle k\rangle^{-M}$. Taking $M>d+\sigma_1-1$, we know that
\begin{align*}
{\rm II}^{(1)}_{m,n}&\lesssim\sum_{n_1,n_2,l\in\mathbb{Z}^d}\langle m-l\rangle^{-N}\bigg(\sum_{k\in\mathbb{Z}^d}\langle k\rangle^{-M}\langle k\rangle^{\sigma_1-1}\bigg)\langle n_1\rangle^{-K}\|\widetilde{\Box}_{l,n_2}u\|_r \chi_{|n-n_1-n_2|_\infty \leq 3 }\nonumber\\
&\lesssim\sum_{n_1,n_2,l\in\mathbb{Z}^d}\langle m-l\rangle^{-N}\langle n_1\rangle^{-K}\|\widetilde{\Box}_{l,n_2}u\|_r \chi_{|n-n_1-n_2|_\infty \leq 3}.
\end{align*}
For $\langle k\rangle \leq \langle l\rangle/10$, we have $\langle l-k\rangle^{-M}\lesssim \langle l\rangle^{-M}$. Then for $M>d+\sigma_1-1$, we can get
\begin{align*}
{\rm II}^{(2)}_{m,n}&\lesssim\sum_{n_1,n_2,l\in\mathbb{Z}^d}\langle m-l\rangle^{-N}\langle l\rangle^{-M}\bigg(\sum_{k\in\mathbb{Z}^d}\langle k\rangle^{\sigma_1-1}\chi_{\{\langle k\rangle \leq \frac{\langle l\rangle}{10}\}}\bigg)\langle n_1\rangle^{-K}\|\widetilde{\Box}_{l,n_2}u\|_r \chi_{|n-n_1-n_2|_\infty \leq 3 }\nonumber\\
&\lesssim\sum_{n_1,n_2,l\in\mathbb{Z}^d}\langle m-l\rangle^{-N}\langle l\rangle^{-M}\langle l\rangle^{d+\sigma_1-1}\langle n_1\rangle^{-K}\|\widetilde{\Box}_{l,n_2}u\|_r \chi_{|n-n_1-n_2|_\infty \leq 3 }\nonumber\\
&\lesssim\sum_{n_1,n_2,l\in\mathbb{Z}^d}\langle m-l\rangle^{-N}\langle n_1\rangle^{-K}\|\widetilde{\Box}_{l,n_2}u\|_r \chi_{|n-n_1-n_2|_\infty \leq 3 }.
\end{align*}
For $\langle k\rangle \in( \langle l\rangle/10, 10\langle l\rangle)$, we need to further divide ${\rm II}^{(3)}_{m,n}$ into two parts,
\begin{align*}
{\rm II}^{(3)}_{m,n}&\lesssim \sum_{n_1,n_2,k,l\in\mathbb{Z}^d}\langle m-l\rangle^{-N}\langle l-k\rangle^{-M}\langle k\rangle^{\sigma_1-1}\langle n_1\rangle^{-K}\|\widetilde{\Box}_{l,n_2}u\|_r \chi_{|n-n_1-n_2|_\infty \leq 3 }\chi_{\{\langle k\rangle > \langle m\rangle\}}\nonumber\\
&\quad\quad\quad\quad\quad \times \chi_{\{\langle k\rangle \in(\frac{\langle l\rangle}{10}, 10\langle l\rangle)\}}\Big(\chi_{\{\langle m\rangle \leq \frac{\langle l\rangle}{10}\}} +\chi_{\{\langle m\rangle \in(\frac{\langle l\rangle}{10}, 10\langle l\rangle)\}}\Big)\nonumber\\
&=:{\rm II}^{(3.1)}_{m,n}+{\rm II}^{(3.2)}_{m,n}.
\end{align*}
If $\langle m\rangle \leq \langle l\rangle/10$, then $\langle m-l\rangle^{-N/2}\lesssim \langle l\rangle^{-N/2}$. For $N/2>d+\sigma_1-1$, we get
\begin{align*}
{\rm II}^{(3.1)}_{m,n}& \lesssim \sum_{n_1,n_2,l\in\mathbb{Z}^d}\langle m-l\rangle^{-N/2}\langle l\rangle^{-N/2}\langle l\rangle^{d+\sigma_1-1}\langle n_1\rangle^{-K}\|\widetilde{\Box}_{l,n_2}u\|_r \chi_{|n-n_1-n_2|_\infty \leq 3 }\nonumber\\
&\lesssim \sum_{n_1,n_2,l\in\mathbb{Z}^d}\langle m-l\rangle^{-N/2}\langle n_1\rangle^{-K}\|\widetilde{\Box}_{l,n_2}u\|_r \chi_{|n-n_1-n_2|_\infty \leq 3 }.
\end{align*}
In the condition of ${\rm II}^{(3.2)}_{m,n}$, we know that $\langle k\rangle\sim \langle l\rangle \sim \langle m\rangle$. Thus for $M>d$,
\begin{align*}
{\rm II}^{(3.2)}_{m,n}& \lesssim \langle m\rangle^{\sigma_1-1}\sum_{n_1,n_2,l\in\mathbb{Z}^d}\langle m-l\rangle^{-N}\bigg(\sum_{k\in\mathbb{Z}^d}\langle l-k\rangle^{-M}\bigg)\langle n_1\rangle^{-K}\|\widetilde{\Box}_{l,n_2}u\|_r \chi_{|n-n_1-n_2|_\infty \leq 3 }\nonumber\\
&\lesssim\langle m\rangle^{\sigma_1-1}\sum_{n_1,n_2,l\in\mathbb{Z}^d}\langle m-l\rangle^{-N}\langle n_1\rangle^{-K}\|\widetilde{\Box}_{l,n_2}u\|_r \chi_{|n-n_1-n_2|_\infty \leq 3 }.
\end{align*}
Now we complete the proof of the Lemma \ref{lemma4}. $\hfill\Box$ \\

{\bf Proof of Lemma \ref{locallemma2}.}  We only prove the result of  \eqref{locallemma2a}, the proof of \eqref{locallemma2aaa} is similar and
easier.   Using Young's inequality, we see that
\begin{align} \label{locallemma2c}
&\|\mathscr{A}_{m,n}((a(x)-a(m))u)\|_{L^\gamma(0,T;L^r)\cap L^\infty(0,T;L^r)}\nonumber\\
&\lesssim \big(\|e^{-tA(m,n)}\|_{L^1(0,T)\cap L^{\gamma'}(0,T)}\big)\cdot\|\Box_{m,n}((a(x)-a(m))u)\|_{L^\gamma(0,T;L^r)}.
\end{align}
Inserting Lemma \ref{lemma4} into \eqref{locallemma2c}, we can get
\begin{align*}
&\|\mathscr{A}_{m,n}((a(x)-a(m))u)\|_{L^\gamma(0,T;L^r)\cap L^\infty(0,T;L^r)}\nonumber\\
&\lesssim \bigg(\int^T_0 \frac{\langle m\rangle^{(\sigma_1-1)\vee0}}{e^{ct|m|^{\sigma_1}}}dt+\bigg(\int^T_0 \frac{\langle m\rangle^{\gamma'(\sigma_1-1)\vee0}}{e^{ct|m|^{\sigma_1}\gamma'}}dt\bigg)^{1/\gamma'}\bigg)\nonumber\\
&\ \ \cdot \sum_{n_1,n_2,l\in\mathbb{Z}^d}\langle m-l\rangle^{-N}\langle n_1\rangle^{-K}\|\widetilde{\Box}_{l,n_2}u\|_{L^\gamma(0,T;L^r)} \chi_{|n-n_1-n_2|_\infty \leq 3 }.
\end{align*}
Applying Minkowski's and Young's inequalities, for $N>d$ we have
\begin{align} \label{locallemma2e}
&\big\|\|\mathscr{A}_{m,n}((a(x)-a(m))u)\|_{L^\gamma(0,T;L^r)\cap L^\infty(0,T;L^r)}\big\|_{\ell^p_m}\nonumber\\
&\lesssim \bigg(\int^T_0 \sup_{m\in\mathbb{Z}^d}\frac{\langle m\rangle^{(\sigma_1-1)\vee0}}{e^{ct|m|^{\sigma_1}}}dt+\bigg(\int^T_0 \sup_{m\in\mathbb{Z}^d}\frac{\langle m\rangle^{\gamma'(\sigma_1-1)\vee0}}{e^{ct|m|^{\sigma_1}\gamma'}}dt\bigg)^{1/\gamma'}\bigg)\nonumber\\
&\quad\quad \cdot\sum_{n_1,n_2\in\mathbb{Z}^d}\langle n_1\rangle^{-K}\Big\|\sum_{l\in\mathbb{Z}^d}\langle m-l\rangle^{-N}\|\widetilde{\Box}_{l,n_2}u\|_{L^\gamma(0,T;L^r)}\Big\|_{\ell^p_m} \chi_{|n-n_1-n_2|_\infty \leq 3 }\nonumber\\
&\lesssim T^\delta \sum_{n_1,n_2\in\mathbb{Z}^d}\langle n_1\rangle^{-K} \big\|\|\Box_{m,n_2}u\|_{L^\gamma(0,T;L^r)}\big\|_{\ell^p_m}\chi_{|n-n_1-n_2|_\infty \leq 3 },
\end{align}
where $T^\delta=\max(T,T^{1/\sigma_1},T^{1/\gamma'},T^{1/\sigma_1-1/\gamma})$. From $\gamma>1\vee\sigma_1$, we see $\delta>0$.
Since $\langle n\rangle^{s_\gamma} \lesssim \langle n_1\rangle^{s_\gamma} + \langle n_2\rangle^{s_\gamma}$ if $|n-n_1-n_2| \lesssim 1$, taking $K>d+s_\gamma$ and using Young's inequality,  we have from \eqref{locallemma2e} that
\begin{align}\label{locallemma2f}
&\|\mathscr{A}_{m,n}((a(x)-a(m))u)\|_{\mathscr{L}^\gamma(0,T;\ell^q_{s_\gamma} \ell^p(L^r))\cap \mathscr{L}^\infty(0,T;\ell^q_{s_0} \ell^p(L^r))}\nonumber\\
&\lesssim T^\delta \bigg(\sum_{n\in\mathbb{Z}^d}\bigg(\sum_{n_1,n_2\in\mathbb{Z}^d}\langle n_1\rangle^{-K+{s_\gamma}} \langle n_2\rangle^{s_\gamma}\big\|\|\Box_{m,n_2}u\|_{L^\gamma(0,T;L^r)}\big\|_{\ell^p_m}\chi_{|n-n_1-n_2|_\infty \leq 3}\bigg)^q\bigg)^{1/q}\nonumber\\
&\lesssim T^\delta\|u\|_{\mathscr{L}^\gamma(0,T;X^{s_\gamma}_{r,p,q})}.
\end{align}
The conclusion follows. $\hfill\Box$

\begin{lem}\label{locallemma3} Let $s\in\mathbb{R}$, $1\leq p,q,r,\gamma\leq\infty$. Assume that $b(\xi)$ satisfies (H1), (H3) and (H4). Then
\begin{gather}
 \|\mathscr{A}_{m,n}((b(D)-b(n))u)\|_{\mathscr{L}^\gamma(0,T;\ell^q_s \ell^p(L^r))}\lesssim \big(T+T^{1/\sigma_2}\big)\|u\|_{\mathscr{L}^\gamma(0,T;X^{s}_{r,p,q})};\label{locallemma3a}\\
 \|\mathscr{A}_{m,n}((b(D)-b(n))u)\|_{L^\gamma(0,T;\ell^q_s \ell^p(L^r))}\lesssim \big(T+T^{1/\sigma_2}\big)\|u\|_{L^\gamma(0,T;X^{s}_{r,p,q})}.\label{locallemma3b}
\end{gather}
\end{lem}
Before proving Lemma \ref{locallemma3}, we show the following lemma.
\begin{lem}\label{lemma5} Let $1\leq r\leq\infty$ and $b(\xi)$ satisfy (H1), (H3) and (H4). Then there exists $\vartheta\in (0\vee(1-\sigma_2),1)$ such that
\begin{align}\label{lemma5a}
 \|\Box_{m,n}((b(D)-b(n))u)\|_r\lesssim \langle n\rangle^{\sigma_2-1}\sum_{|l_1|_\infty\leq1}\sum_{k\in\mathbb{Z}^d}\langle m-k\rangle^{-(d+1-\vartheta)} \|\Box_{k,n+l_1}u\|_r.
\end{align}
\end{lem}
{\bf Proof.}  Noticing  the almost orthogonality $\sigma_n=\sum_{|l_1|_\infty\leq1}\sigma_n\sigma_{n+l_1}$, we see that
\begin{align}\label{lemma5d}
&\|\Box_{m,n}((b(D)-b(n))u)\|_r\nonumber\\
&\lesssim \sum_{|l_1|_\infty \vee |l_2|_\infty\leq1}\sum_{k\in \mathbb{Z}^d}\|\sigma_m(x) \mathscr{F}^{-1} \sigma_{n}(\xi)(b(\xi)-b(n)) \mathscr{F}\sigma_{k+l_2}(x)\Box_{k,n+l_1}u\|_r.
\end{align}
When $n\neq0$, from Lemma \ref{lemma1} and \eqref{remark1}, for $N\in\mathbb{N}$ and $\theta\in[0,1]$, we obtain
\begin{align}\label{lemma5b}
&\|\Box_{m,n}((b(D)-b(n))u)\|_r\nonumber\\
&\lesssim\sum_{|l_1|_\infty\leq1}\sum_{k\in\mathbb{Z}^d}\|\sigma\|_r\|\sigma\|_{r'}\sup_{|\beta|\leq N}\big\|\partial^\beta_\xi\big(\sigma_{n}(\xi)(b(\xi)-b(n))\big)\big\|_1\langle m-k\rangle^{-\theta N} \|\Box_{k,n+l_1}u\|_r.
\end{align}
From the condition (UD), the mean value theorem and the assumption of $b(\xi)$, it follows that
\begin{align*}
\big\|\partial^\beta_\xi\big(\sigma_{n}(\xi)(b(\xi)-b(n))\big)\big\|_\infty\lesssim \langle n\rangle^{\sigma_2-1},\ \   \forall \ |\beta|\leq N.
\end{align*}
Since $\sigma_{n}(\xi)$ has compact support, we have from \eqref{lemma5b} that
\begin{align}\label{lemma5c}
\|\Box_{m,n}((b(D)-b(n))u)\|_r
\lesssim \langle n\rangle^{\sigma_2-1}\sum_{|l_1|_\infty\leq1}\sum_{k\in\mathbb{Z}^d}\langle m-k\rangle^{-\theta N} \|\Box_{k,n+l_1}u\|_r.
\end{align}
When $n=0$, similar  to the proof of Lemma \ref{lemma1}, we need to control $${\rm I}(x,y):=\bigg|\int\sigma(\xi)(b(\xi)-b(0))e^{{\rm i}(x-y)\cdot \xi}d\xi\bigg|.$$
If $|x-y|_\infty\lesssim 1$, ${\rm I}(x,y)\leq \|\sigma(\xi)(b(\xi)-b(0))\|_1\lesssim 1$ because of $b(\xi)\sim |\xi|^{\sigma_2}$, $|\xi|\leq 1$. If $|x-y|_\infty\gg1$, we suppose $|x_1-y_1|=|x-y|_\infty\gg1$. For the sake of convenience, denote $\bar{b}(\xi):=b(\xi)-b(0)$ and  $\bar{\xi}=(\xi_2, \xi_3, \cdots, \xi_d)$. Using dyadic decomposition, we have
\begin{align}\label{frac0}
{\rm I}(x,y)&\leq\sum_{j\leq0}\bigg|\int\sigma(\xi)\varphi_j(\xi)\bar{b}(\xi)e^{{\rm i}(x-y)\cdot \xi}d\xi\bigg|\nonumber\\
&=\sum_{j\leq0}2^{jd}\bigg|\int_{\frac{1}{2}\leq|\xi|\leq2}\sigma(2^j\xi)\varphi(\xi)\bar{b}(2^j\xi)e^{{\rm i}2^j(x_1-y_1)\cdot \xi_1}d\xi_1e^{{\rm i}2^j(\bar{x}-\bar{y})\cdot \bar{\xi}}d\bar{\xi}\bigg|\nonumber\\
&=: \sum_{j\leq0}2^{jd}\cdot{\rm J}(x,y).
\end{align}
By carrying out the repeated integrations by parts, we know that
\begin{align*}
{\rm J}(x,y)
&\lesssim\frac{1}{|2^j(x_1-y_1)|^K}\int_{\frac{1}{2}\leq|\xi|\leq2}\bigg|\frac{\partial^K(\sigma(2^j\xi)
\varphi(\xi)\bar{b}(2^j\xi))}{\partial\xi_1^K}\bigg|d\xi_1d\bar{\xi}\nonumber\\
&\lesssim 2^{j(\sigma_2-K)}|x_1-y_1|^{-K}, \ \ \ \ \forall K\in\mathbb{N}^+.
\end{align*}
Hence, making the interpolation between $K=d$ and $K=d+1$, we have for any $\vartheta\in[0,1)$,
\begin{align}\label{frac1}
{\rm J}(x,y)&\lesssim 2^{j(\sigma_2-d)\vartheta}|x_1-y_1|^{-d\vartheta}2^{j(\sigma_2-(d+1))(1-\vartheta)}|x_1-y_1|^{-(d+1)(1-\vartheta)}\nonumber\\
&\lesssim 2^{j(\sigma_2-(d+1-\vartheta))}|x_1-y_1|^{-(d+1-\vartheta)}.
\end{align}
Inserting \eqref{frac1} into \eqref{frac0} and taking $1>\vartheta>1-\sigma_2$, we get
\begin{align*}
{\rm I}(x,y)\lesssim \sum_{j\leq0}2^{j(\sigma_2-1+\vartheta)}\cdot|x_1-y_1|^{-(d+1-\vartheta)}\lesssim |x_1-y_1|^{-(d+1-\vartheta)}\lesssim \langle x-y\rangle^{-(d+1-\vartheta)}.
\end{align*}
Therefore, \eqref{lemma5d} is continued by
\begin{align}\label{Estimate10}
\|\Box_{m,0}((b(D)-b(0))u)\|_r\lesssim \sum_{|l_1|_\infty\leq1}\sum_{k\in\mathbb{Z}^d}\langle m-k \rangle ^{-(d+1-\vartheta)} \|\Box_{k,l_1}u\|_r.
\end{align}
Combining \eqref{lemma5c} and \eqref{Estimate10}, taking $\theta N=d+1-\vartheta$, we get the result \eqref{lemma5a}. $\hfill\Box$\\

{\bf Proof of Lemma \ref{locallemma3}.} We only prove the result of \eqref{locallemma3a}, because \eqref{locallemma3b} is similar. From Lemma \ref{lemma5} and using Young's inequality,  we have
\begin{align*}
&\|\mathscr{A}_{m,n}((b(D)-b(n))u)\|_{L^\gamma(0,T;L^r)}\nonumber\\
&\lesssim
\int^T_0e^{-t\mathfrak{Re}A(m,n)}dt\cdot\sum_{|l_1|_\infty\leq1}\sum_{k\in\mathbb{Z}^d}\langle n\rangle^{\sigma_2-1}\langle m-k\rangle^{-(d+1-\vartheta)}\|\Box_{k,n+l_1}u\|_{L^\gamma(0,T;L^r)}
\end{align*}
By Young's inequality and $d+1-\vartheta>d$, we obtain
\begin{align}\label{Estimate3}
&\Big\|\|\mathscr{A}_{m,n}((b(D)-b(n))u)\|_{L^\gamma(0,T;L^r)}\Big\|_{\ell_m^p}\nonumber\\
&\lesssim \bigg(\int^T_0\frac{\langle n\rangle^{\sigma_2-1}}{e^{ct| n|^{\sigma_2}}}dt\bigg)
\Big\|\sum_{|l_1|_\infty\leq1}\sum_{k\in\mathbb{Z}^d}\langle m-k\rangle^{-(d+1-\vartheta)}\|\Box_{k,n+l_1}u\|_{L^\gamma(0,T;L^r)}\Big\|_{\ell_m^p}\nonumber\\
&\lesssim \bigg(\int^T_0\frac{\langle n\rangle^{\sigma_2-1}}{e^{ct| n|^{\sigma_2}}}dt\bigg)\Big\|\sum_{|l_1|_\infty\leq1}\|\Box_{m,n+l_1}u\|_{L^\gamma(0,T;L^r)}\Big\|_{\ell_m^p}
\end{align}
From the definition we can obtain from \eqref{Estimate3} that
\begin{align*}
&\|\mathscr{A}_{m,n}((b(D)-b(n))u)\|_{\mathscr{L}^\gamma(0,T;\ell^q_s \ell^p(L^r))}\nonumber\\
&\lesssim \bigg(\int^T_0\sup_{n\in\mathbb{Z}^d}\frac{\langle n\rangle^{\sigma_2-1}}{e^{ct| n|^{\sigma_2}}}dt\bigg)\bigg\|\langle n\rangle^{s}\Big\|\sum_{|l_1|_\infty\leq1}\|\Box_{m,n+l_1}u\|_{L^\gamma(0,T;L^r)}\Big\|_{\ell_m^p}\bigg\|_{\ell_n^q}\nonumber\\
&\lesssim \big(T+T^{1/\sigma_2}\big)\|u\|_{\mathscr{L}^\gamma(0,T;X^{s}_{r,p,q})}.
\end{align*}
The conclusion follows. $\hfill\Box$

\section{Nonlinear Estimates}

\begin{prop}\label{nonl}
Let   $s\in\mathbb{R}$, $1\leq p,q,r,\gamma,\gamma_1,\gamma_2\leq\infty$,  $1+1/\gamma=1/\gamma_1+1/\gamma_2$, and $0\leq\kappa\leq\sigma_2/\gamma_2$, $0<T<1$. Then we have
\begin{align} \label{nonlinear7}
\|\mathscr{A}_{m,n}f\|& _{\mathscr{L}^\gamma(0,T;\ell^q_s \ell^p (L^r)) \cap \mathscr{L}^\infty(0,T;\ell^q_{s-\sigma_2/\gamma} \ell^p (L^r)) } \nonumber\\
&  \lesssim  (T^{1/\gamma'_1} + T^{1/\gamma_2-\kappa/\sigma_2}) \|f\|_{\mathscr{L}^{\gamma_1}(0,T;X^{s-\kappa}_{r,p,q})}.
\end{align}

\end{prop}
{\bf Proof.} Using Young's inequality, we have
\begin{gather*}
 \|\mathscr{A}_{m,n}f\|_{L^\gamma(0,T;L^r)}\lesssim \|e^{-tA(m,n)}\|_{L^{\gamma_2}_{t\in[0,T]}}\|\Box_{m,n}f\|_{L^{\gamma_1}(0,T;L^r)};\\
 \langle n\rangle^{-\sigma_2/\gamma}\|\mathscr{A}_{m,n}f\|_{L^\infty(0,T;L^r)}\lesssim \langle n\rangle^{-\sigma_2/\gamma} \|e^{-tA(m,n)}\|_{L^{\gamma'_1}_{t\in[0,T]}}\|\Box_{m,n}f\|_{L^{\gamma_1}(0,T;L^r)}.
 \end{gather*}
Since  $\mathfrak{Re}A(m,n)\gtrsim |n|^{\sigma_2}$, we see that
\begin{gather*}
\|e^{-tA(m,n)}\|_{L^{\gamma_2}_{t\in[0,T]}}\lesssim \|e^{-ct|n|^{\sigma_2}}\|_{L^{\gamma_2}_{t\in[0,T]}}\lesssim \langle n\rangle^{-\sigma_2/\gamma_2};\\
\langle n\rangle^{-\sigma_2/\gamma}\|e^{-tA(m,n)}\|_{L^{\gamma'_1}_{t\in[0,T]}}\lesssim\langle n\rangle^{-\sigma_2(1/\gamma+1/\gamma'_1)} \lesssim \langle n\rangle^{-\sigma_2/\gamma_2}.
\end{gather*}
On the other hand, if $0\leq\kappa<\sigma_2/\gamma_2$, we have
\begin{align*}
\|e^{-tA(m,n)}\|_{L^{\gamma_2}_{t\in[0,T]}}&\lesssim \langle n\rangle^{-\kappa}\bigg(\int_0^T \frac{\langle n\rangle^{\kappa\gamma_2}}{ e^{ct|n|^{\sigma_2}\gamma_2}}dt\bigg)^{1/{\gamma_2}}\\
&\lesssim \langle n\rangle^{-\kappa}\bigg(\int_0^T 1+ t^{-\frac{\kappa\gamma_2}{\sigma_2}}dt\bigg)^{1/{\gamma_2}}
\lesssim T^{1/\gamma_2-\kappa/\sigma_2}\langle n\rangle^{-\kappa},
\end{align*}
and
\begin{align*}
\langle n\rangle^{-\sigma_2/\gamma}\|e^{-tA(m,n)}\|_{L^{\gamma'_1}_{t\in[0,T]}}&\lesssim \langle n\rangle^{-\kappa}\bigg(\int_0^T \Big(\frac{\langle n\rangle^{\kappa-\sigma_2/\gamma}}{ e^{ct|n|^{\sigma_2}}}\Big)^{\gamma'_1}dt\bigg)^{1/{\gamma'_1}}\\
&\lesssim \langle n\rangle^{-\kappa}\bigg(\int_0^T 1+ t^{-\big(\frac{\kappa}{\sigma_2}-\frac{1}{\gamma}\big)\gamma'_1}dt\bigg)^{1/{\gamma'_1}} \nonumber\\
& \lesssim ( T^{1/\gamma'_1}  + T^{1/\gamma_2-\kappa/\sigma_2}) \langle n\rangle^{-\kappa}.
\end{align*}
Hence, we have
\begin{align} \label{nonlinear6}
\|\mathscr{A}_{m,n}f\|_{L^\gamma(0,T;L^r)} & +\langle n\rangle^{-\sigma_2/\gamma}\|\mathscr{A}_{m,n}f\|_{L^\infty(0,T;L^r)} \nonumber\\
& \lesssim  ( T^{1/\gamma'_1}  + T^{1/\gamma_2-\kappa/\sigma_2}) \langle n\rangle^{-\kappa}\|\Box_{m,n}f\|_{L^{\gamma_1}(0,T;L^r)}.
\end{align}
Taking sequence $\ell^q_s\ell^p$  norms in both sides of \eqref{nonlinear6}, we get \eqref{nonlinear7}, as desired. $\hfill\Box$\\

\begin{lem}\label{derivative} {\rm (Derivative Estimates)}
Let $s\in\mathbb{R}$, $1\leq p,q,r\leq\infty$. We have
\begin{gather}
\|\partial_x^\alpha u\|_{\mathscr{L}^\gamma(0,T; X^s_{r,p,q})}\lesssim \|u\|_{\mathscr{L}^\gamma(0,T; X^{s+|\alpha|}_{r,p,q})}.\label{Estimate17}
\end{gather}
\end{lem}
{\bf Proof.} It is easy to see that
\begin{align}
&\big\|\Box_{m,n}(\partial_x^\alpha u)\big\|_r=\|\sigma_m(x)\mathscr{F}^{-1}\sigma_n(\xi)\xi^\alpha \hat{u}\|_r\nonumber\\
&\lesssim\sum_{|l_1|_\infty\vee|l_2|_\infty\leq1}\sum_{k\in\mathbb{Z}^d}\|\sigma_m(x)\mathscr{F}^{-1}\sigma_{n+l_1}(\xi)\xi^\alpha \mathscr{F}\sigma_{k+l_2}\sigma_k \mathscr{F}^{-1}\sigma_n \hat{u}\|_r. \label{Estimate11}
\end{align}
From Lemma \ref{lemma1}, \eqref{Estimate11} is controlled by
\begin{align} \label{multiplier}
\big\|\Box_{m,n}(\partial_x^\alpha u)\big\|_r &\lesssim\sum_{|l_1|_\infty\leq1}\sum_{k\in\mathbb{Z}^d}\|\sigma\|_r\|\sigma\|_{r'}\sup_{|\beta|\leq N}\|\partial^\beta_\xi(\sigma_{n+l_1}(\xi)\xi^\alpha)\|_1\langle m-k\rangle^{-N}\|\Box_{k,n}u\|_r\nonumber\\
&\lesssim \langle n\rangle^{|\alpha|}\sum_{k\in\mathbb{Z}^d}\langle m-k\rangle^{-N}\|\Box_{k,n}u\|_r.
\end{align}
 Inserting \eqref{multiplier} into the definition of $X^s_{r,p,q}$, using Young's inequality and taking $N>d$, we can obtain
\begin{align*}
\|\partial_x^\alpha u\|_{ \mathscr{L}^\gamma(0,T;X^s_{r,p,q})} &=\Bigg(\sum_{ n\in\mathbb{Z}^d} \langle n\rangle^{sq}\bigg(\sum_{ m\in\mathbb{Z}^d} \big\|\Box_{m,n}(\partial_x^\alpha u)\big\|^p_{L^\gamma(0,T;L^r)}\bigg)^{q/p}\Bigg)^{1/q}\nonumber\\
&\lesssim \Bigg(\sum_{ n\in\mathbb{Z}^d} \langle n\rangle^{(s+|\alpha|)q} \bigg(\sum_{ m\in\mathbb{Z}^d} \bigg(\sum_{k\in\mathbb{Z}^d}\langle m-k\rangle^{-N}\|\Box_{k,n}u\|_{L^\gamma(0,T;L^r)}\bigg)^p\bigg)^{q/p}\Bigg)^{1/q}\nonumber\\
&\lesssim \|u\|_{ \mathscr{L}^\gamma(0,T;X^{s+|\alpha|}_{r,p,q})} .
\end{align*}
Now we get the conclusion \eqref{Estimate17}.  $\hfill\Box$

\begin{lem}\label{locallemma4} Let $1\leq p,q,r,\gamma\leq\infty$, $\kappa:=\max\limits_{1\leq i\leq N}|\alpha_i|\leq \sigma_2$  and $ 0<T<1$.
\begin{itemize}
  \item [\rm (i)]  If $\kappa<\sigma_2$, let  $s_\gamma>\kappa+d/q'$,   $\gamma > \max(N,(N-1)\sigma_2/(\sigma_2-\kappa))$, then there exists $\delta(N)= \min(  1/\gamma'_1, \ 1+(1-N)/\gamma -\kappa/\sigma_2 )>0$ such that
      $$\|\mathscr{A}_{m,n}(\partial_x^{\alpha_1} v_1 \cdots \partial_x^{\alpha_N} v_N )\|_{\mathscr{L}^\gamma(0,T;\ell^q_{s_\gamma} \ell^p(L^r))\cap\mathscr{L}^\infty(0,T;\ell^q_{s_0}\ell^p(L^r))}\lesssim T^{\delta(N)} \prod^N_{i=1} \|v_i\|_{\mathscr{L}^\gamma(0,T;X^{s_\gamma}_{r,p,q})}.$$
  \item [\rm (ii)] If $\kappa=\sigma_2$, let $s_0>\sigma_2+d/q'$,  then
  $$\|\mathscr{A}_{m,n}(\partial_x^{\alpha_1} v_1 \cdots \partial_x^{\alpha_N} v_N )\|_{\mathscr{L}^\infty(0,T;\ell^q_{s_0}\ell^p(L^r))}\lesssim \prod^N_{i=1}\|v_i\|_{\mathscr{L}^\infty(0,T;X^{s_0}_{r,p,q})}.$$
  \end{itemize}
\end{lem}
 {\bf Proof.} Denote $u_1=\partial_x^{\alpha_1} v_1, \cdots, u_N=\partial_x^{\alpha_N} v_N$. From Propositon \ref{nonl}, we have
\begin{align}
&\|\mathscr{A}_{m,n}(u_1 u_2\cdots u_N)\|_{\mathscr{L}^\gamma(0,T;\ell^q_{s_\gamma} \ell^p(L^r))\cap\mathscr{L}^\infty(0,T;\ell^q_{s_0} \ell^p(L^r))}\nonumber\\
&\lesssim T^{\delta(N)}  \|u_1u_2\cdots u_N\|_{\mathscr{L}^{\gamma_1}(0,T;X^{s_\gamma-\kappa}_{r,p,q})}.
\end{align}
where $1+1/\gamma=1/\gamma_1+1/\gamma_2$, $0\leq\kappa\leq\sigma_2/\gamma_2$.
Then from Proposition \ref{nonlinear}, Lemma \ref{derivative}, we have
\begin{align}\label{nonlinear5}
&\|\mathscr{A}_{m,n}(u_1 u_2\cdots u_N)\|_{\mathscr{L}^\gamma(0,T;\ell^q_{s_\gamma} \ell^p(L^r))\cap\mathscr{L}^\infty(0,T;\ell^q_{s_0} \ell^p(L^r))}\nonumber\\
&\lesssim T^{\delta(N)} \prod^N_{i=1}\|u_i\|_{\mathscr{L}^\gamma(0,T;X^{s_\gamma-\kappa}_{r,p,q})}\lesssim  T^{\delta(N)} \prod^N_{i=1} \|v_i\|_{\mathscr{L}^\gamma(0,T;X^{s_\gamma}_{r,p,q})},
\end{align}
where $1/\gamma_1=N/\gamma$, $\gamma> N$. If $\kappa<\sigma_2$,  from $\gamma>N \vee (N-1)\sigma_2/(\sigma_2-\kappa)$
it follows that
\begin{align*}
  \frac{1}{\gamma_2}-\frac{\kappa}{\sigma_2}>0, \ \  \frac{1}{\gamma'_1}  >0.
\end{align*}
Then we obtain the conclusion of (i). If $\kappa=\sigma_2$, it means that $\gamma_2=1$, $\gamma=\gamma_1=\infty$. Note that $s_\infty=s_0$, then \eqref{nonlinear5} implies the conclusion (ii).  $\hfill\Box$\\

In order to handle the case that the nonlinearity contains no derivative, we need a modification to the above argument. Let us connect our approach with Proposition \ref{nonl}. We consider the following two cases:

{\it Case 1.} $ 1/q \geq \sigma_2/d $. Let $s>0$ satisfy
$$
s>d/q' -\sigma_2/(N-1).
$$
Then we can choose some $\varepsilon_i \in (0,1)$ satisfying
$$
\frac{s(1-\varepsilon_0)}{d} =\frac{1}{q'} -\frac{\sigma_2(1-\varepsilon_1) (1-\varepsilon_2)}{d(N-1)}
$$
Denote $\kappa=\sigma_2 (1-\varepsilon_2)$, and taking $\gamma >(N-1)/ \varepsilon_2$ (one can assume that $\varepsilon_2<(N-1)/N$),  $\gamma_1 =\gamma/N$,  we have from Proposition \ref{nonl} that
\begin{align} \label{nonlinear7+}
\|\mathscr{A}_{m,n}f\| _{\mathscr{L}^\gamma(0,T;\ell^q_s \ell^p (L^r)) \cap \mathscr{L}^\infty(0,T;\ell^q_{s-\sigma_2/\gamma} \ell^p (L^r)) }
  \lesssim  T^{\delta}  \|f\|_{\mathscr{L}^{\gamma/N}(0,T; X^{s-\kappa}_{r,p,q})}
\end{align}
for some $\delta>0.$ We can choose $q_1$ and $q_2$ satisfying
$$
\frac{1}{q} - \frac{1}{q_1} = \frac{\sigma_2(1-\varepsilon_1) (1-\varepsilon_2)}{d}, \ \  \frac{1}{q'_2} = \frac{\sigma_2(1-\varepsilon_1) (1-\varepsilon_2)}{d(N-1)}.
$$
Applying H\"older's inequality, we have from \eqref{nonlinear7+} that
\begin{align} \label{nonlinear7++}
\|\mathscr{A}_{m,n}f\| _{\mathscr{L}^\gamma(0,T;\ell^q_s \ell^p (L^r)) \cap \mathscr{L}^\infty(0,T;\ell^q_{s-\sigma_2/\gamma} \ell^p (L^r)) }
  \lesssim  T^{\delta}   \|f\|_{\mathscr{L}^{\gamma/N}(0,T; X^{s}_{r,p,q_1})}.
\end{align}
Noticing that
$$
\frac{1}{q_1}  = \frac{1}{q }  +  \frac{N-1}{q_2}   -(N-1) ,
$$
in view of Proposition \ref{nonlinear}, we have
\begin{align} \label{nonlinear7+3}
   \|u^N \|_{\mathscr{L}^{\gamma/N}(0,T; X^{s}_{r,p,q_1})}  \lesssim \|u\|_{\mathscr{L}^\gamma (0,T; X^s_{r,p,q})} \|u\|^{N-1}_{\mathscr{L}^\gamma (0,T; X^0_{\infty,\infty,q_2})}.
\end{align}
Since $(1/q_2-1/q)=s(1-\varepsilon_0)/d$, we have the embedding
$$
\mathscr{L}^\gamma (0,T; X^s_{r, p, q}) \subset \mathscr{L}^\gamma (0,T; X^0_{\infty,\infty,q_2}).
$$
Hence,
\begin{align} \label{nonlinear7+4}
 \|u^N \|_{\mathscr{L}^{\gamma/N}(0,T; X^{s}_{r,p,q_1})}   \lesssim \|u\|^N_{\mathscr{L}^\gamma (0,T; X^s_{r,p,q})}.
\end{align}
It follows from \eqref{nonlinear7++} and  \eqref{nonlinear7+4} that
\begin{align} \label{nonlinear7+5}
\|\mathscr{A}_{m,n}(u^N)\| _{\mathscr{L}^\gamma(0,T;\ell^q_s \ell^p (L^r)) \cap \mathscr{L}^\infty(0,T;\ell^q_{s-\sigma_2/\gamma} \ell^p (L^r)) }
     \lesssim  T^\delta  \|u\|^N_{\mathscr{L}^\gamma (0,T; X^s_{r,p,q})}.
\end{align}

{\it Case 2.} We consider the case $1/q< \sigma_2/d$.  Let $s>0$ satisfy
$$
s>d/q' - d/q(N-1).
$$
We can find some $\varepsilon_i \in (0,1)$ satisfying
$$
\frac{s(1-\varepsilon_0)}{d} =\frac{1}{q'} -\frac{(1-\varepsilon_1)}{ q (N-1)}.
$$
By choosing
\begin{align} \label{gamma}
\gamma > \max\left(N, \ \ \frac{N-1}{1-d/q\sigma_2} \right),
\end{align}
we have from Proposition \ref{nonl} that
\begin{align} \label{nonlinear7+a}
\|\mathscr{A}_{m,n}f\| _{\mathscr{L}^\gamma(0,T;\ell^q_s \ell^p (L^r)) \cap \mathscr{L}^\infty(0,T;\ell^q_{s-\sigma_2/\gamma} \ell^p (L^r)) }
  \lesssim  T^{\delta}   \|f\|_{\mathscr{L}^{\gamma/N}(0,T; X^{s-d/q}_{r,p,q})}
\end{align}
for some $\delta>0.$ We can choose $q_1$ and $q_2$ satisfying
$$
\frac{1}{q} - \frac{1}{q_1} = \frac{ (1-\varepsilon_1)  }{q}, \ \  \frac{1}{q'_2} = \frac{ (1-\varepsilon_1)}{q(N-1)}.
$$
Applying H\"older's inequality, we have from \eqref{nonlinear7+} that
\begin{align} \label{nonlinear7++a}
\|\mathscr{A}_{m,n}f\| _{\mathscr{L}^\gamma(0,T;\ell^q_s \ell^p (L^r)) \cap \mathscr{L}^\infty(0,T;\ell^q_{s-\sigma_2/\gamma} \ell^p (L^r)) }
  \lesssim  T^{\delta}  \|f\|_{\mathscr{L}^{\gamma/N}(0,T; X^{s}_{r,p,q_1})}.
\end{align}
Noticing that
$$
\frac{1}{q_1}  = \frac{1}{q }  +  \frac{N-1}{q_2}   -(N-1) ,
$$
in view of Proposition \ref{nonlinear}, we have
\begin{align} \label{nonlinear7+3a}
 \|u^N \|_{\mathscr{L}^{\gamma/N}(0,T; X^{s}_{r,p,q_1})}  \lesssim \|u\|_{\mathscr{L}^\gamma (0,T; X^s_{r,p,q})} \|u\|^{N-1}_{\mathscr{L}^\gamma (0,T; X^0_{\infty,\infty,q_2})}.
\end{align}
Since $(1/q_2-1/q)=s(1-\varepsilon_0)/d$, we have the embedding
$
\mathscr{L}^\gamma (0,T; X^s_{r, p, q}) \subset \mathscr{L}^\gamma (0,T; X^0_{\infty,\infty,q_2}).
$
Hence,
it follows that \eqref{nonlinear7+5} also holds in Case 2.  Up to now, we have shown the following:

\begin{lem}  \label{powerest}
Assume that $s\geq 0$ satisfies
$$
s>
\left\{
\begin{array}{ll}
d/q' -\sigma_2/(N-1), & \ if \ \sigma_2\leq d/q,\\
d/q' - d/q(N-1), & \ \ if \ \sigma_2>d/q.
\end{array}
\right.
$$
Let $\gamma \gg N$ if $\sigma_2\leq d/q$, and $\gamma$ satisfy \eqref{gamma} if $\sigma_2> d/q$. Then there exists $\delta>0$ such that
 \begin{align} \label{nonlinear7+5a}
\|\mathscr{A}_{m,n}(v_1...v_N)\| _{\mathscr{L}^\gamma(0,T;\ell^q_s \ell^p (L^r)) \cap \mathscr{L}^\infty(0,T;\ell^q_{s-\sigma_2/\gamma} \ell^p (L^r)) }
     \lesssim T^{\delta} \prod^N_{i=1} \|v_i\|_{\mathscr{L}^\gamma (0,T; X^s_{r,p,q})}.
\end{align}
\end{lem}

\section{Local Well-Posedness for $A(x,D)=a(x)+b(D)$}

{\bf Proof of (i) and (ii) of Theorem \ref{thma} (Local Well-Posedness).}  The conclusion is obtained by using the standard iteration method. Let us take $F(u)= \lambda_1 \partial_x^{\alpha_1} u\cdots \partial_x^{\alpha_K} u +  \lambda_2 \partial_x^{\beta_1} u  \cdots \partial_x^{\beta_{L}} u$ and the general cases can be handled in a similar way. We can assume that $K\geq L$.
We consider the iteration sequence $\{u^{(\mu)}\}^\infty_{\mu=0}$:
\begin{align}
\partial_t u^{(\mu+1)} +A(x,D) u^{(\mu+1)} =F(u^{(\mu)}), \ \ u^{(0)}=0.
\end{align}
Using Proposition \ref{discrit integral},
\begin{align}\label{Estimate22}
 \Box_{m,n}  u^{(\mu+1)} (t) = &   e^{-t A(m,n)}\Box_{m,n} u_0-   \mathscr{A}_{m,n}((a(x)-a(m))u^{(\mu+1)})\nonumber\\
  &-   \mathscr{A}_{m,n}((b(D)-b(n))u^{(\mu+1)})+   \mathscr{A}_{m,n}(F(u^{(\mu)})), \ \ \ u^{(0)}=0.
\end{align}
We will show that $\{u^{(\mu)}\}^\infty_{\mu=0}$ is a Cauchy sequence in  $\mathscr{L}^\gamma (0,T; X^{s_\gamma}_{r,p,q})$ for some $T>0$.

{\rm (i)} $0<\kappa<\sigma_2$. Recall that $\gamma$ satisfies
$$
 \gamma > K \vee  \frac{(K-1)\sigma_2}{\sigma_2-\kappa} \geq L \vee  \frac{(L-1)\sigma_2}{\sigma_2-\kappa} .
$$
Assume that $s_\gamma=s_0+\sigma_2/\gamma>\kappa+d/q'$.  We can assume that $0<T<1$. Taking $N=K,L$ in Lemma \ref{locallemma4}, we have
\begin{align}
\|\mathscr{A}_{m,n} & (F(v)  -  F(w))\|_{\mathscr{L}^\gamma(0,T;\ell^q_{s_\gamma} \ell^p(L^r))} \nonumber\\
& \lesssim T^{\delta(L)} (\|v\|^{L-1}_{\mathscr{L}^\gamma(0,T;X^{s_\gamma}_{r,p,q})} + \|w\|^{L-1}_{\mathscr{L}^\gamma(0,T;X^{s_\gamma}_{r,p,q})})\|v -w\|_{\mathscr{L}^\gamma(0,T; X^{s_\gamma}_{r,p,q})} \nonumber\\
 & \quad  +  T^{\delta(K)} (\|v\|^{K-1}_{\mathscr{L}^\gamma(0,T;X^{s_\gamma}_{r,p,q})} +  \|w\|^{K-1}_{\mathscr{L}^\gamma(0,T;X^{s_\gamma}_{r,p,q})}  ) \|v -w\|_{\mathscr{L}^\gamma(0,T; X^{s_\gamma}_{r,p,q})}. \label{limit1}
\end{align}
Combining \eqref{Estimate22}, Lemmas \ref{locallemma2}, \ref{locallemma3} and the above estimate,  we obtain that there exists $\delta>0$ such that
\begin{align}
\|\mathscr{A}_{m,n} (A(x,D) -A(m,n))(v  &-  w )\|_{\mathscr{L}^\gamma(0,T;\ell^q_{s_\gamma} \ell^p(L^r))} \nonumber\\
\lesssim& T^\delta  \|v - w\|_{\mathscr{L}^\gamma(0,T; X^{s_\gamma}_{r,p,q})}. \label{limit2}
\end{align}
Taking $T>0$ satisfies $CT^\delta \leq 1/2$, we have from \eqref{limit1} and \eqref{limit2} that
\begin{align}
\|& u^{(\mu+1)} - u^{(\mu)}\|_{\mathscr{L}^\gamma(0,T; X^{s_\gamma}_{r,p,q})} \nonumber\\
& \lesssim T^{\delta(L)} (\|u^{(\mu)}\|^{L-1}_{\mathscr{L}^\gamma(0,T;X^{s_\gamma}_{r,p,q})} + \|u^{(\mu-1)}\|^{L-1}_{\mathscr{L}^\gamma(0,T;X^{s_\gamma}_{r,p,q})})\|u^{(\mu)}-u^{(\mu-1)}\|_{\mathscr{L}^\gamma(0,T; X^{s_\gamma}_{r,p,q})} \nonumber\\
 & \quad  +  T^{\delta(K)} (\|u^{(\mu)}\|^{K-1}_{\mathscr{L}^\gamma(0,T;X^{s_\gamma}_{r,p,q})} +  \|u^{(\mu-1)}\|^{K-1}_{\mathscr{L}^\gamma(0,T;X^{s_\gamma}_{r,p,q})}  ) \|u^{(\mu)}-u^{(\mu-1)}\|_{\mathscr{L}^\gamma(0,T; X^{s_\gamma}_{r,p,q})}. \label{limit3}
\end{align}
Take $M=3C\|u_0\|_{M^{s_0}_{p,q}}$. We further assume that $C (T^{\delta(L)} M^{L-1} + T^{\delta(K)} M^{K-1}) \leq 1/100 $. It follows from Lemma \ref{locallemma1} that
$$
\|u^{(1)} \|_{\mathscr{L}^\gamma(0,T; X^{s_\gamma}_{r,p,q})} \leq C \|u_0\|_{M^{s_0}_{p,q}} = M/3.
$$
By induction, we have
$$
\|u^{(\mu)} \|_{\mathscr{L}^\gamma(0,T; X^{s_\gamma}_{r,p,q})} \leq   M
$$
and
\begin{align*}
\|u^{(\mu+1)} & - u^{(\mu)}\|_{\mathscr{L}^\gamma(0,T; X^{s_\gamma}_{r,p,q})} \leq \frac{1}{2} \|u^{(\mu)}-u^{(\mu-1)}\|_{\mathscr{L}^\gamma(0,T; X^{s_\gamma}_{r,p,q})}.
\end{align*}
This implies that $\{u^{(\mu)}\} $ is a Cauchy sequence in $\mathscr{L}^\gamma(0,T; X^{s_\gamma}_{r,p,q})$ and must converge to some $u \in \mathscr{L}^\gamma(0,T; X^{s_\gamma}_{r,p,q})$. Taking $w=u$ and $v=u^{(\mu)}$ in \eqref{limit1};  and taking $w=u$ and $v=u^{(\mu+1)}$ in \eqref{limit2}, we immediately have from \eqref{Estimate22} that
\begin{align}\label{limit22}
 \Box_{m,n}  u (t) = &   e^{-t A(m,n)}\Box_{m,n} u_0-   \mathscr{A}_{m,n}((a(x)-a(m))u )\nonumber\\
  &-   \mathscr{A}_{m,n}((b(D)-b(n))u )+   \mathscr{A}_{m,n} F(u )
\end{align}
for all $m,n\in \mathbb{Z}^d$.

 We can further show that $L^\infty(0,T; X^{s_0}_{r,p,q})$.    Indeed, from Minkowski's inequality, and Lemmas \ref{locallemma1}, \ref{locallemma2}, \ref{locallemma3} and \ref{locallemma4}, we see that
\begin{align*}
 \|u\|_{L^\infty(0,T; X^{s_0}_{r,p,q})}
&\leq C\|u_0\|_{M^{s_0}_{p,q}}+CT^\delta \|u\|_{L^\infty(0,T; X^{s_0}_{r,p,q})}+ CT^{\delta}\|u\|_{\mathscr{L}^\gamma(0,T;X^{s_\gamma}_{r,p,q})} \nonumber\\
&  \quad +CT^{\delta}\|u\|^K_{\mathscr{L}^\gamma(0,T;X^{s_\gamma}_{r,p,q})}+ CT^{\delta} \|u\|^L_{\mathscr{L}^\gamma(0,T; X^{s_\gamma}_{r,p,q})}.
\end{align*}
From the above discussion, we have
\begin{align*}
\|u\|_{L^\infty(0,T; X^{s_0}_{r,p,q})} \lesssim   \|u_0\|_{M^{s_0}_{p,q}}.
\end{align*}
Therefore, the system \eqref{limit22} has a unique solution $u \in C(0,T; X^{s_0}_{r,p,q})\cap \mathscr{L}^\gamma (0,T; X^{s_0+\sigma_2/\gamma}_{r,p,q})$.

{\rm (ii)} $\kappa=\sigma_2$. Let $\gamma=\infty$ and $s_0>\sigma_2+d/q'$. Take $M=3C\|u_0\|_{ M^{s_0}_{p,q}} $. Let $T$ and $\|u_0\|_{M^{s_0}_{p,q}}$ be sufficiently small such that
$$
CT^\delta + C(M^{K-1} +M^{L-1})\leq 1/100,
$$
  Similarly as above,
combining \eqref{Estimate22}, Lemmas \ref{locallemma1}, \ref{locallemma2}, \ref{locallemma3} and \ref{locallemma4}, we obtain
\begin{align*}
\|u^{(\mu+1)} & - u^{(\mu)}\|_{\mathscr{L}^\infty(0,T; X^{s_0}_{r,p,q})} \nonumber\\
 & \lesssim    (\|u^{(\mu)} \|^{K-1}_{\mathscr{L}^\infty(0,T;X^{s_0}_{r,p,q})} + \|u^{(\mu-1)} \|^{K-1}_{\mathscr{L}^\infty(0,T;X^{s_0}_{r,p,q})}) \|u^{(\mu)} - u^{(\mu-1)} \|_{\mathscr{L}^\infty(0,T;X^{s_0}_{r,p,q})} \nonumber\\
& \quad  +  (\|u^{(\mu)} \|^{L-1}_{\mathscr{L}^\infty(0,T;X^{s_0}_{r,p,q})} + \|u^{(\mu-1)} \|^{L-1}_{\mathscr{L}^\infty(0,T;X^{s_0}_{r,p,q})}) \|u^{(\mu)} - u^{(\mu-1)} \|_{\mathscr{L}^\infty(0,T;X^{s_0}_{r,p,q})}.
\end{align*}

  Repeating the arguments as in the case $\kappa<\sigma_2$, we can get that there exists a unique solution $u$ in $C(0,T; X^{s_0}_{r,p,q})\cap \mathscr{L}^\infty (0,T; X^{s_0}_{r,p,q})$. $\hfill\Box$\\

{\bf Proof of (i) and (ii) of Theorem \ref{thma0} (Local Well-Posedness).} By Lemma \ref{powerest}, one can use the same way  as in the proof of the local well-posedness of Theorem \ref{thma} to obtain the result and the details of the proof are omitted.  $\hfill\Box$

\section{Global well-posedness for $A(x,D)=a(x)+b(D)$}

\subsection{Initial data in $L^2$}
We need a priori estimate.
\begin{lem}\label{priori}{\rm (A Priori Estimate)} Assume that $A(x,D)$ satisfies $\mathfrak{Re}(A(x,D)u,u)\geq \|u\|^2_{\dot{H}^{\sigma_2/2}}$,  and $F(u)$ satisfies $\mathfrak{Re}(F(u),u)\leq 0$. Then we have a priori estimate
\begin{align}\label{priori0}
\int_0^t\|u(\tau)\|^2_{\dot{H}^{\sigma_2/2}}d\tau+ \frac{1}{2} \|u(t)\|^2_2\leq \frac{1}{2} \|u_0\|^2_2.
\end{align}
\end{lem}
{\bf Proof.} Taking the inner-product with $u$ in \eqref{PDE}, we get
\begin{align}\label{global1}
\big(\partial_t u  + A(x,D)u -F(u), u\big)=0.
\end{align}
Taking the real part of \eqref{global1} and using the assumptions, we have
\begin{align*}
0=\mathfrak{ Re}\big(\partial_t u,u\big)  +\mathfrak{ Re} \big(A(x,D)u,u\big) - \mathfrak{ Re}\big( F(u), u\big)
\geq \frac{1}{2} \partial_t \|u\|^2_2+\|u\|^2_{\dot{H}^{\sigma_2/2}}.
\end{align*}
Integrating about the time variable, we know that \eqref{priori0} holds. $\hfill\Box$\\

{\bf Proof of (iii) of Theorem \ref{thma0}.} First, let us observe the priori estimate.  From condition (H4) it follows that,
\begin{align*}
{\mathfrak{Re}}(A(x,D)u,u)&= {\mathfrak{Re}} \int a(x)|u|^2 dx+  {\mathfrak{Re}}\int b(\xi)|\hat{u}|^2 d\xi\nonumber\\
& \gtrsim \int |x|^{\sigma_1}|u|^2 dx + \int |\xi|^{\sigma_2}|\hat{u}|^2 d\xi\geq\|u\|^2_{\dot{H}^{\sigma_2/2}}.
\end{align*}
It is easy to see that  $\mathfrak{Re}(\lambda_i|u|^{k_i-1}u, \ u)\leq 0$. Thus from the  priori estimate \eqref{priori0}, we see that $\sup_{t\in[0,\infty]}\|u(t)\|_2\leq\|u_0\|_2$.

 Recall that $M^0_{2,2} =L^2$ with equivalent norm.
If $\sigma_2> d/2$, $s> d(1-1/(K-1))/2$, using the result of (ii) of Theorem \ref{thma0}, we can take
$$
\gamma=
\left\{
\begin{array}{ll}
K+ \varepsilon, & K\leq 2\sigma_2/d,\\
(K-1)/(1-d/2\sigma_2)+ \varepsilon, & K>  2\sigma_2/d
\end{array}\right.
$$
for some $0<\varepsilon \ll 1.$ If $K<1+ 2\sigma_2/d$, then we have
$$
\frac{d}{2}\left(1- \frac{1}{K-1}\right) -\frac{\sigma_2}{\gamma} <0.
$$
It follows that there exists $s>d(1- 1/(K-1))/2 $ such that $s-\sigma_2/\gamma =0$.  So, in view of the result of (ii) of Theorem \ref{thma0}, $u_0\in L^2$ implied that there exists a $T>0$ such that \eqref{PDE} has a unique solution $u\in C([0,T]; L^2)\cap \mathscr{L}^\gamma (0,T; X^{\sigma_2/\gamma }_{2, 2,2})$.  Since $\sup_{t\in[0,\infty)} \|u(t)\|_2\leq\|u_0\|_2$, we see that one can extend the solution from $[0,T]$ to $[T, 2T]$, $[2T, 3T]$,... .
$\hfill\Box$

\subsection{Derivative Nonlinearity }

Generally speaking, if $\mathfrak{Re}a(0)= \mathfrak{Re}b(0)=0$, $\partial_t + A(m,n)$ contains no dissipative structure at $(m,n)=(0,0)$. As a result, it is impossible to get a time-global estimate for $\mathscr{A}_{0,0} f$. Hence, we assume $A(x,\xi)$ satisfies (H1)-(H4) and in addition, assume $a(0)>0$ in this subsection.   Recall that $s_\gamma=s_0+\sigma_2/\gamma$, $s_\infty=s_0$ and $s_2=s_0+\sigma_2/2$.

\begin{lem}\label{globallemma1} Let $1\leq p,q,r\leq\infty$ and $\gamma\in \{2,\infty\}$.  We have
\begin{align*}
   \|e^{-t A(m,n)}\Box_{m,n} u_0\|_{\mathscr{L}^\gamma(0,\infty;\ell^q_{s_\gamma} \ell^p(L^r))}\lesssim \|u_0\|_{X^{s_0}_{r,p,q}}.
\end{align*}
\end{lem}
{\bf Proof.}  It is easy to show that
\begin{align*}
\|e^{-t A(m,n)}\Box_{m,n} u_0\|_{L^\gamma(0,\infty;L^r)}&\leq\bigg(\int_0^\infty e^{-t{ \mathfrak{Re}A(m,n)}\gamma } dt\bigg)^{1/\gamma} \|\Box_{m,n} u_0\|_{L^r}\nonumber\\
&\lesssim \langle n\rangle^{-\sigma_2/\gamma}\|\Box_{m,n} u_0\|_{L^r}.
\end{align*}
By the definition \eqref{defn}, we get the desired conclusion. $\hfill\Box$

\begin{lem}\label{globallemma2} Let  $1\leq p,q,r\leq\infty$, $s_0\geq 0$. Assume that $a(x)$ satisfies (H1), (H2), (H4) and $a(0)>0 $. Then there exists $\rho>0$ such that
\begin{gather}
\|\mathscr{A}_{m,n}((a(x)-a(m))u)\|_{\mathscr{L}^2(0,\infty;\ell^q_{s_2} \ell^p(L^r))}\lesssim \|u\|_{\mathscr{L}^2(0,\infty;X^{s_2-\rho}_{r,p,q})};\label{globallemma2a}\\
\|\mathscr{A}_{m,n}((a(x)-a(m))u)\|_{\mathscr{L}^\infty(0,\infty;\ell^q_{s_0} \ell^p(L^r))}
\lesssim \|u\|_{\mathscr{L}^{2\vee\sigma_1}(0,\infty;X^{s_0}_{r,p,q})}.\label{globallemma2b}
\end{gather}
\end{lem}
{\bf Proof.} For \eqref{globallemma2a},
using Young's inequality, we know that
\begin{align}\label{globallemma2c}
&\|\mathscr{A}_{m,n}((a(x)-a(m))u)\|_{L^2_tL_x^r}\nonumber\\
&\lesssim \int^\infty_0 e^{-t \mathfrak{Re}A(m,n)} dt\cdot\|\Box_{m,n}((a(x)-a(m))u)\|_{L^2_tL_x^r}\nonumber\\
&\lesssim \frac{1}{\langle m\rangle^{\sigma_1}+\langle n\rangle^{\sigma_2}}\cdot\|\Box_{m,n}((a(x)-a(m))u)\|_{L^2_tL_x^r},
\end{align}
where the last inequality is because of $a(0)\geq c$ and $ \mathfrak{Re}A(x,\xi)\gtrsim   |x|^{\sigma_1}+|\xi|^{\sigma_2}$. In view of Young's inequality $ab\leq a^p/p+b^{p'}/p'$ ($a,b>0$, $p>1$), we see that $\langle m\rangle^{\sigma_1}+\langle n\rangle^{\sigma_2} \gtrsim (\langle m\rangle^{\sigma_1})^{(\sigma_1-1)/\sigma_1}(\langle n\rangle^{\sigma_2})^{1/\sigma_1}\gtrsim \langle m\rangle^{\sigma_1-1} \langle n\rangle^{\sigma_2/\sigma_1} $, if $\sigma_1>1$.  Therefore, we have the following facts:
\begin{align*}
\frac{1}{\langle m\rangle^{\sigma_1}+\langle n\rangle^{\sigma_2} }\lesssim
\left\{
    \begin{array}{ll}
      \langle n\rangle^{-\sigma_2} , & \hbox{if $0\leq\sigma_1\leq 1$;}\vspace{1ex} \\
      \langle m\rangle^{1-\sigma_1} \langle n\rangle^{-\sigma_2/\sigma_1}, & \hbox{if $\sigma_1>1$.}
    \end{array}
  \right.
\end{align*}
Therefore, \eqref{globallemma2c} implies that
\begin{align}\label{globallemma2d}
&\|\mathscr{A}_{m,n}((a(x)-a(m))u)\|_{L^2_tL_x^r}\nonumber\\
&\lesssim \frac{1}{\langle m\rangle^{(\sigma_1-1)\vee0}}\langle n\rangle^{-\sigma_2/1\vee\sigma_1}\cdot\|\Box_{m,n}((a(x)-a(m))u)\|_{L^2_tL_x^r},
\end{align}
Inserting \eqref{lemma4a} into \eqref{globallemma2d}, we can get
\begin{align*}
&\|\mathscr{A}_{m,n}((a(x)-a(m))u)\|_{L^2_tL_x^r}\nonumber\\
&\lesssim \langle n\rangle^{-\sigma_2/1\vee\sigma_1}\sum_{n_1,n_2,l\in\mathbb{Z}^d}\langle m-l\rangle^{-N}\langle n_1\rangle^{-K}\|\widetilde{\Box}_{l,n_2}u\|_{L^2_tL_x^r} \chi_{|n-n_1-n_2|\leq k_0(d)}.
\end{align*}
By similar discussion with \eqref{locallemma2e}-\eqref{locallemma2f}, we immediately get the conclusion \eqref{globallemma2a}.

For \eqref{globallemma2b}, if $\sigma_1\leq 2$, then $\sigma_1/2\geq (\sigma_1-1)\vee0$, from Young's inequality we get
\begin{align}\label{globallemma2f}
\|\mathscr{A}_{m,n}((a(x)-a(m))u)\|_{L^\infty_tL_x^r}&\lesssim \langle m\rangle^{-\sigma_1/2}
\|\Box_{m,n}((a(x)-a(m))u)\|_{L^2_tL_x^r}\nonumber\\
&\lesssim \langle m\rangle^{-((\sigma_1-1)\vee0)}
\|\Box_{m,n}((a(x)-a(m))u)\|_{L^2_tL_x^r}.
\end{align}
If $\sigma_1>2$, then from Young's inequality and $1+1/\infty=1/\sigma_1+(\sigma_1-1)/\sigma_1$, we know that
\begin{align}\label{globallemma2g}
\|\mathscr{A}_{m,n}((a(x)-a(m))u)\|_{L^\infty_tL_x^r}
\lesssim \langle m\rangle^{-(\sigma_1-1)}
\|\Box_{m,n}((a(x)-a(m))u)\|_{L^{\sigma_1}_tL_x^r}.
\end{align}
Combining \eqref{globallemma2f}-\eqref{globallemma2g} and Lemma \ref{lemma4}, we have
\begin{align*}
&\|\mathscr{A}_{m,n}((a(x)-a(m))u)\|_{L^\infty_tL_x^r}\nonumber\\
&\lesssim \sum_{n_1,n_2,l\in\mathbb{Z}^d}\langle m-l\rangle^{-N}\langle n_1\rangle^{-K}\|\widetilde{\Box}_{l,n_2}u\|_{L^{2\vee\sigma_1}_tL_x^r} \chi_{|n-n_1-n_2|\leq k_0(d)}.
\end{align*}
Again by similar discussion with \eqref{locallemma2e}-\eqref{locallemma2f}, we have the conclusion \eqref{globallemma2b}.

\begin{lem}\label{globallemma3} Let $1\leq p,q,r\leq\infty$ and $\gamma\in \{2,\infty\}$. Assume that $b(\xi)$ satisfies (H1), (H3) and (H4). Then
\begin{align}
 \|\mathscr{A}_{m,n}((b(D)-b(n))u)\|_{\mathscr{L}^\gamma(0,\infty;\ell^q_{s_\gamma} \ell^p(L^r))}\lesssim \|u\|_{\mathscr{L}^2(0,\infty;X^{s_2-1}_{r,p,q})}.
\end{align}
\end{lem}
{\bf Proof.} By using Young's inequality and $1+1/\gamma=1/2+(\gamma+2)/2\gamma$, we have
\begin{align}\label{globallemma3a}
&\|\mathscr{A}_{m,n}((b(D)-b(n))u)\|_{L^\gamma_tL^r_x}\nonumber\\
&\lesssim\big\|e^{-tA(m,n)}\|_{L^{2\gamma/\gamma+2}_t}\|\Box_{m,n}((b(D)-b(n))u)\|_{L^2_tL^r_x}\nonumber\\
&\lesssim\langle n\rangle^{-\sigma_2(\gamma+2)/2\gamma}\|\Box_{m,n}((b(D)-b(n))u)\|_{L^2_tL^r_x}.
\end{align}
Note that $s_\gamma-\sigma_2(\gamma+2)/2\gamma+\sigma_2-1=s_2-1$ for $\gamma\in \{2,\infty\}$, then using Lemma \ref{lemma5} and Young's inequality, we have from \eqref{globallemma3a} that
\begin{align*}
&\|\mathscr{A}_{m,n}((b(D)-b(n))u)\|_{\mathscr{L}^\gamma(0,\infty;\ell^q_{s_\gamma} \ell^p(L^r))}\nonumber\\
&\lesssim \bigg\|\langle n\rangle^{s_2-1}\Big\|\sum_{|l_1|_\infty\leq1}\sum_{k\in\mathbb{Z}^d}\langle m-k\rangle^{-(d+1-\vartheta)}\|\Box_{k,n+l_1}u\|_{L^\gamma_tL^r_x}\Big\|_{\ell_m^p}\bigg\|_{\ell_n^q}\nonumber\\
&\lesssim \bigg\|\langle n\rangle^{s_2-1}\Big\|\sum_{|l_1|_\infty\leq1}\|\Box_{m,n+l_1}u\|_{L^\gamma_tL^r_x}\Big\|_{\ell_m^p}\bigg\|_{\ell_n^q}
\lesssim \|u\|_{\mathscr{L}^\gamma(0,\infty;X^{s_2-1}_{r,p,q})},
\end{align*}
where we used the condition $d+1-\vartheta>d$. The conclusion follows. $\hfill\Box$
\begin{lem}\label{globallemma4} Let $1\leq p,q,r\leq\infty$, $\kappa:=\max\limits_{1\leq i\leq N}|\alpha_i|\leq\sigma_2$, and $s_0>\kappa+d/q'$, then we have
\begin{align*}
&\|\mathscr{A}_{m,n}(\partial_x^{\alpha_1} u\cdot \partial_x^{\alpha_2} u\cdots \partial_x^{\alpha_N} u )\|_{\mathscr{L}^2(0,\infty;\ell^q_{s_2}\ell^p(L^r))\cap\mathscr{L}^\infty(0,\infty;\ell^q_{s_0}\ell^p(L^r))}\nonumber\\
&\lesssim \|u\|^{N-1}_{\mathscr{L}^\infty(0,\infty;X^{s_0}_{r,p,q})}\|u\|_{\mathscr{L}^2(0,\infty;X^{s_2}_{r,p,q})}.
\end{align*}
\end{lem}
{\bf Proof.} Denote $u_1=\partial_x^{\alpha_1} u$,  $\cdots$, $u_N=\partial_x^{\alpha_N} u$. From Young's inequality, we have
\begin{gather*}
 \|\mathscr{A}_{m,n}(u_1 u_2\cdots u_N)\|_{L^2_tL^r_x}\lesssim
 \langle n\rangle^{-\sigma_2}\|\Box_{m,n}(u_1u_2\cdots u_N)\|_{L^{2}_tL^r_x};\\
 \|\mathscr{A}_{m,n}(u_1 u_2\cdots u_N)\|_{L^\infty_tL^r_x}\lesssim
 \langle n\rangle^{-\sigma_2/2}\|\Box_{m,n}(u_1u_2\cdots u_N)\|_{L^{2}_tL^r_x}.
\end{gather*}
It follows from Proposition \ref{nonlinear} and Lemma \ref{derivative} that
\begin{align*}
&\|\mathscr{A}_{m,n}(u_1 u_2\cdots u_N)\|_{\mathscr{L}^2(0,\infty;\ell^q_{s_2}\ell^p(L^r))\cap\mathscr{L}^\infty(0,\infty;\ell^q_{s_0}\ell^p(L^r))}\nonumber\\
&\lesssim \Bigg(\sum_{n\in\mathbb{R}^d}\langle n\rangle^{(s_2-\sigma_2)q}\bigg(\sum_{m\in\mathbb{R}^d}\|\Box_{m,n}(u_1u_2\cdots u_N)\|^p_{L^{2}_tL^r_x}\bigg)^{q/p}\Bigg)^{1/q}\nonumber\\
& \lesssim \sum^N_{j=1}\prod_{i\neq j}\|u_i\|_{\mathscr{L}^\infty(0,\infty;X^{d/q'+}_{r,p,q})}
\|u_j\|_{\mathscr{L}^2(0,\infty;X^{s_2-\kappa}_{r,p,q})}  \nonumber\\
&\lesssim \|u\|^{N-1}_{\mathscr{L}^\infty(0,\infty;X^{s_0}_{r,p,q})}\|u\|_{\mathscr{L}^2(0,\infty;X^{s_2}_{r,p,q})},
\end{align*}
The conclusion follows. $\hfill\Box$\\

{\bf Proof of (iii) of Theorem \ref{thma}.} Without loss of generality we can assume that $F(u)=\lambda_1 \partial_x^{\alpha_1} u \cdots \partial_x^{\alpha_K} u + \lambda_2 \partial_x^{\alpha_1} u\cdots \partial_x^{\alpha_L} u$. Let us assume that $\|u_0\|_{H^{s_0}} \leq \delta$ for some sufficieintly small $\delta>0$. In view of $X^s_{r,2,2}=H^s$, we will denote $\mathscr{L}^q(0,T; H^{s}):= \mathscr{L}^q(0,T; X^{s}_{r,2,2})$. According to the local well-posedness results, it suffices to show that
\begin{align} \label{bound}
\varphi(T):=  \|u\|_{\mathscr{L}^\infty(0,T; H^{s_0}) \cap \mathscr{L}^2(0,T; H^{s_2})}\leq C \delta, \ \ \forall \ T>0.
\end{align}
Recall the integral equation
\begin{align*}
 \Box_{m,n} u (t) = & e^{-t A(m,n)}\Box_{m,n} u_0-\mathscr{A}_{m,n}((a(x)-a(m))u)\nonumber\\
  &-\mathscr{A}_{m,n}((b(D)-b(n))u)+\mathscr{A}_{m,n}(F(u)).
\end{align*}
By Lemma \ref{globallemma1}--Lemma \ref{globallemma4}, we obtain
\begin{align}\label{2Estimate4}
&\|u\|_{\mathscr{L}^\infty(0,T; H^{s_0})}+\|u\|_{\mathscr{L}^2(0,T; H^{s_2})}\nonumber\\
&\lesssim \|u_0\|_{H^{s_0}}+\|u\|_{\mathscr{L}^2(0,T; H^{s_2-\rho})}+\|u\|_{\mathscr{L}^{2\vee\sigma_1}(0,T; H^{s_0})}\nonumber\\
&\quad +(\|u\|^{K-1}_{\mathscr{L}^\infty(0,T;H^{s_0})} + \|u\|^{L-1}_{\mathscr{L}^\infty(0,T;H^{s_0})}) \|u\|_{\mathscr{L}^2(0,T;H^{s_2})}.
\end{align}
 Since $a(0)>0$ and $ \mathfrak{Re} A(x,\xi) \gtrsim \langle x \rangle^{\sigma_1}+|\xi|^{\sigma_2}$, we have $\mathfrak{Re}(A(x,D)u,u)\geq \|u\|^2_{{H}^{\sigma_2/2}}$. In a similar way to a priori estimate \eqref{priori0}, we can get
\begin{align*}
\int_0^t\|u(\tau)\|^2_{{H}^{\sigma_2/2}}d\tau+ \frac{1}{2} \|u(t)\|^2_2\leq \frac{1}{2} \|u_0\|^2_2.
\end{align*}
It follows that
$$\|u\|_{\mathscr{L}^2(0,T;H^{\sigma_2/2})}=\|u\|_{L^2(0,T;H^{\sigma_2/2})}\lesssim \|u_0\|_2.\nonumber$$
\indent If $\sigma_1\leq 2$,  by interpolation theory, we have from \eqref{2Estimate4} that
\begin{align}\label{2Estimate1}
&\|u\|_{\mathscr{L}^\infty(0,T; H^{s_0})}+\|u\|_{\mathscr{L}^2(0,T; H^{s_2})}\nonumber\\
&\lesssim \|u_0\|_{H^{s_0}}+\|u\|_{\mathscr{L}^2(0,T;H^{s_2-\rho})}+(\|u\|^{K-1}_{\mathscr{L}^\infty(0,T;H^{s_0})} + \|u\|^{L-1}_{\mathscr{L}^\infty(0,T;H^{s_0})}) \|u\|_{\mathscr{L}^2(0,T;H^{s_2})}\nonumber\\
&\lesssim \|u_0\|_{H^{s_0}}+\|u\|^\theta_{\mathscr{L}^2(0,T;H^{\sigma_2/2})}\|u\|^{1-\theta}_{\mathscr{L}^2(0,T;H^{s_2})}
+ \varphi(T)^{K}  + \varphi(T)^{L}\nonumber\\
&\lesssim \|u_0\|_{H^{s_0}}+\|u_0\|^\theta_{2} \varphi(T)^{1-\theta} + \varphi(T)^{K}  + \varphi(T)^{L}.
\end{align}
The continuity arguments imply that \eqref{bound} holds. Indeed, assume that $\|u_0\|_{H^{s_0}} \leq \epsilon$ for some $0< \epsilon\ll 1$.  If $\varphi(T) \leq \epsilon^\theta$, then from \eqref{2Estimate1} it follows that
\begin{align}\label{2Estimatecont}
\varphi(T)
 \lesssim \epsilon +\epsilon^\theta  \epsilon^{\theta(1-\theta)} + \epsilon^{\theta K}  + \epsilon^{\theta L}.
\end{align}
Since $\epsilon\ll 1$, we have $\varphi(T) \leq \epsilon^\theta/2$. This implies that $\varphi(T) \leq \epsilon^\theta $  for all $T>0$.

 If $\sigma_1 >2$, we have
\begin{align*}
\|u\|_{\mathscr{L}^{\sigma_1}(0,T;H^{s_0})}&\lesssim \|u\|^{2/\sigma_1}_{\mathscr{L}^2(0,T;H^{s_0})}\|u\|^{(\sigma_1-2)/\sigma_1}_{\mathscr{L}^\infty(0,T;H^{s_0})}\nonumber\\
&\lesssim \|u\|^{2\vartheta/\sigma_1}_{\mathscr{L}^2(0,T;H^{\sigma_2/2})}\|u\|^{2(1-\vartheta)/\sigma_1}_{\mathscr{L}^2(0,T;H^{s_2})}\|u\|^{(\sigma_1-2)/\sigma_1}_{\mathscr{L}^\infty(0,\infty;H^{s_0})}\nonumber\\
&\lesssim \|u_0\|^{2\vartheta/\sigma_1}_{2}\varphi(T)^{1-2\vartheta/\sigma_1}.
\end{align*}
Similar to \eqref{2Estimate1}, we have from \eqref{2Estimate4} that
\begin{align*}
 \| u\|_{\mathscr{L}^\infty(0,T; H^{s_0}) \cap \mathscr{L}^2(0,T; H^{s_2})}
\lesssim & \ \|u_0\|_{H^{s_0}}+\|u_0\|^\theta_{2}\varphi(T)^{1-\theta} \nonumber\\
& +\|u_0\|^{2\vartheta/\sigma_1}_{2}\varphi(T)^{1-2\vartheta/\sigma_1}+\varphi(T)^L + \varphi(T)^K.
\end{align*}
Again, in view of the continuity arguments, one can obtain that \eqref{bound} holds. $\hfill\Box$

\section{Well-posedness in Case   $A\in S^M$ }

In this section we will always assume that
\begin{itemize}
\item[(A1)]  $A\in S^M$  is $M$-order dissipative;

\item[(A2)]  $|\partial_x^\alpha A(x,\xi)|\leq A_\alpha(1+|\xi|)^{M-\varepsilon} \ for \   \ 1\leq |\alpha|\leq N.$
\end{itemize}

We only need to estimate the second term of the localized integral equation in Proposition \ref{discrit integral},   the estimates of the other terms have been worked out in the previous sections.

\begin{lem}\label{lemmaB}  Let $1\leq p,q\leq\infty$, $\partial^\alpha_{x}\partial^\beta_{\xi}\varphi(x,\xi)\in L^1_{loc}(\mathbb{R}^d\times\mathbb{R}^d)$.
Denote for $m,n,h,l,k\in \mathbb{Z}^d$,
\begin{align*}
T_{m,n,h,l,k}f(x):=\Box_{m,n}\big(\sigma_{h}(x)\mathscr{F}^{-1}_\xi\sigma_{l}(\xi)\varphi(x,\xi)\mathscr{F}_x\sigma_{k}(x)f(x)\big).
\end{align*}
Then for any $N_1,N_2,N_3\in\mathbb{N}$, we have
\begin{align}\label{lemmaB0}
\|T_{m,n,h,l,k}f\|_q\lesssim \langle m-h\rangle^{-{N_1}}\langle n-l\rangle^{-{N_2}}\langle h-k\rangle^{-{N_3}}\sup\limits_{|\alpha|\leq N_2\atop |\beta|\leq N_3}\|\partial^\alpha_{x}\partial^\beta_{\xi}\varphi(x,\xi)\|_{L^1_{x,\xi}(Q_h\times Q_l)} \|f\|_p.
\end{align}
\end{lem}
{\bf Proof.} Let us rewrite
\begin{align}\label{lemmaB1}
T_{m,n,h,l,k}f(x)=\sigma_m(x) \iiiint e^{{\rm i}x\eta}\sigma_n(\eta)e^{-{\rm i}z\eta}\sigma_{h}(z)e^{{\rm i}z\xi}\sigma_{l}(\xi)\varphi(z,\xi)e^{-{\rm i}y\xi}\sigma_{k}(y)f(y)dyd\xi dzd\eta.
\end{align}
\emph{Case} 1. We consider the case $|m_1-h_1|=|m-h|_\infty\gg 1$, $|h_1-k_1|=|h-k|_\infty\gg 1$ and $|n_1-l_1|=|n-l|_\infty\gg 1$. Making the integration by part, we have
\begin{align}\label{lemmaB2}
\int e^{{\rm i}(x-z)\eta}\sigma_n(\eta)d\eta=\frac{1}{(-{\rm i})^{N_1}(x_1-z_1)^{N_1}}\int e^{{\rm i}(x-z)\eta}\partial^{N_1}_{\eta_1}\sigma_n(\eta)d\eta,
\end{align}
\begin{align}\label{lemmaB3}
&\int e^{{\rm i}(\xi-\eta)z}(x_1-z_1)^{-{N_1}}\sigma_{h}(z)\varphi(z,\xi)dz\nonumber\\
&=\frac{1}{(-{\rm i})^{N_2}(\xi_1-\eta_1)^{N_2}}\int e^{{\rm i}(\xi-\eta)z}\partial^{N_2}_{z_1}\big((x_1-z_1)^{-{N_1}}\sigma_{h}(z)\varphi(z,\xi)\big)dz,
\end{align}
and
\begin{align}\label{lemmaB4}
&\int e^{{\rm i}(z-y)\xi}(\xi_1-\eta_1)^{-{N_2}}\sigma_{l}(\xi)\partial^{N_2}_{z_1}\big((x_1-z_1)^{-{N_1}}\sigma_{h}(z)\varphi(z,\xi)\big)d\xi\nonumber\\
&=\frac{1}{(-{\rm i})^{N_3}(z_1-y_1)^{N_3}}\int e^{{\rm i}(z-y)\xi}\partial^{N_3}_{\xi_1}\Big((\xi_1-\eta_1)^{-{N_2}}\sigma_{l}(\xi)\partial^{N_2}_{z_1}
\big((x_1-z_1)^{-{N_1}}\sigma_{h}(z)\varphi(z,\xi)\big)\Big)d\xi.
\end{align}
Inserting \eqref{lemmaB2}-\eqref{lemmaB4} into \eqref{lemmaB1}, and using the condition (UD), we obtain
\begin{align*}
|T_{m,n,h,l,k}f(x)|&\leq \bigg|\sigma_m(x) \iiiint|\partial^{N_1}_{\eta_1}\sigma_n(\eta)||z_1-y_1|^{-{N_3}}|\sigma_{k}(y)f(y)|\nonumber\\
&\quad\times\Big|\partial^{N_3}_{\xi_1}\Big((\xi_1-\eta_1)^{-{N_2}}\sigma_{l}(\xi)\partial^{N_2}_{z_1}
\big((x_1-z_1)^{-{N_1}}\sigma_{h}(z)\varphi(z,\xi)\big)\Big)\Big|dyd\xi dzd\eta\bigg|\nonumber\\
&\lesssim |m_1-h_1|^{-{N_1}}|n_1-l_1|^{-{N_2}}|h_1-k_1|^{-{N_3}}|\sigma_m(x)|\int|\partial^{N_1}_{\eta_1}\sigma_n(\eta)|d\eta\nonumber\\
&\quad \times \sup\limits_{|\alpha|\leq N_2\atop |\beta|\leq N_3} \int_{Q_h}\int_{Q_l} |\partial^\beta_{\xi_1}\partial^\alpha_{z_1}\varphi(z,\xi)|d\xi dz\int|\sigma_{k}(y)f(y)|dy.
\end{align*}
Therefore, by using H$\ddot{\rm o}$lder's inequality, the result \eqref{lemmaB0} follows.\\
\emph{Case} 2. At least one of $|m-h|_\infty$, $|h-k|_\infty$ and $|n-l|_\infty$ is less than $C$. Without loss of generality, we assume that $|n_1-l_1|=|n-l|_\infty\lesssim 1$. In this case we do not need to make the integrations by parts of \eqref{lemmaB3}. We can get the corresponding result
\begin{align*}
\|T_{m,n,h,l,k}f\|_q\lesssim \langle m-h\rangle^{-{N_1}}\langle h-k\rangle^{-{N_3}}\sup\limits_{|\beta|\leq N_3}\|\partial^\beta_{\xi}\varphi(x,\xi)\|_{L^1_{x,\xi}(Q_h\times Q_l)} \|f\|_p.
\end{align*}
This completes the proof of Lemma \ref{lemmaB}. $\hfill\Box$

In the following we will use the same notations as in the previous section by letting $s_\gamma:=s_0+M/\gamma$, i.e. $s_\infty=s_0$, $s_2=s_0+M/2$.

\begin{lem}\label{lemmaA}  Let  $1\leq p,q,r,\gamma\leq\infty$, $s_\gamma\geq 0$. Assume that $A(x,\xi)$ satisfies (A1) and (A2), then we have
\begin{align*}
&\|\mathscr{A}_{m,n}((A(x,D)-A(m,n))u)\|_{\mathscr{L}^\gamma(0,T;\ell^q_{s_\gamma} \ell^p(L^r))\cap \mathscr{L}^\infty(0,T;\ell^q_{s_0} \ell^p(L^r))}\nonumber\\
&\lesssim \big(T+T^{\varepsilon/M}+T^{1/\gamma'}\big)\|u\|_{\mathscr{L}^\gamma(0,T;X^{s_\gamma}_{r,p,q})}.
\end{align*}
\end{lem}
{\bf Proof.} By the properties of pseudo-differential operators and $\sum_{l\in\mathbb{Z}^d} \sigma_{l}(x)=1$, we have
\begin{align}\label{lemmaA1}
&\Box_{m,n}((A(x,D)-A(m,n))u)\nonumber\\
=&\ \sigma_m(x) \mathscr{F}_\eta^{-1}\sigma_n(\eta)\mathscr{F}_z\mathscr{F}_\xi^{-1}( A(z,\xi) -A(m,n))\hat{u}\nonumber\\
=&\sum_{h,l\in\mathbb{Z}^d}\sigma_m(x) \mathscr{F}_\eta^{-1}\sigma_n(\eta)\mathscr{F}_z\sigma_{h}(z)\mathscr{F}_\xi^{-1}\sigma_{l}(\xi)( A(z,\xi) -A(m,n))\hat{u}.
\end{align}
Denote
\begin{align*}
\overline{\Box}_{k,l}:=\sum_{|l_1|_\infty\vee|l_2|_\infty\leq1}\sigma_{k+l_1}(x)\mathscr{F}^{-1}\sigma_{l+l_2}(\xi)\mathscr{F},
\end{align*}
we have from \eqref{lemmaA1} that
\begin{align*}
&\Box_{m,n}((A(x,D)-A(m,n))u)\nonumber\\
=&\sum_{h,k,l\in\mathbb{Z}^d}\sigma_m(x) \mathscr{F}_\eta^{-1}\sigma_n(\eta)\mathscr{F}_z\sigma_{h}(z)\mathscr{F}_\xi^{-1}\sigma_{l}(\xi)( A(z,\xi) -A(m,n))\mathscr{F}_y\sigma_{k}(y)\overline{\Box}_{k,l}u.
\end{align*}
From Lemma \ref{lemmaB},
\begin{align}\label{lemmaA2}
&\|\Box_{m,n}((A(x,D)-A(m,n))u)\|_r\nonumber\\
&\lesssim\sum_{h,k,l\in\mathbb{Z}^d}\langle m-h\rangle^{-{N_1}-1}\langle n-l\rangle^{-{N_2}-1}\langle h-k\rangle^{-{N_3}}\nonumber\\
&\qquad\times\sup\limits_{|\alpha|\leq N_2\atop |\beta|\leq N_3}\|\partial^\beta_{\xi}\partial^\alpha_{z}( A(z,\xi) -A(m,n))\|_{L^1_{z,\xi}(Q_h\times Q_l)}\|\overline{\Box}_{k,l}u\|_r.
\end{align}
From the condition (UD), the mean value theorem and the assumption of $A(x,\xi)$, it follows that
\begin{align}\label{lemmaA3}
&\|\partial^\beta_{\xi}\partial^\alpha_{z}( A(z,\xi) -A(m,n))\|_{L^1_{z,\xi}(Q_h\times Q_l)}\nonumber\\
&\lesssim \big(\langle l\rangle^{M-\varepsilon}+\langle n\rangle^{M-\varepsilon}\big)\langle m-h\rangle+\big(\langle l\rangle^{M-1}+\langle n\rangle^{M-1}\big)\langle n-l\rangle.
\end{align}
Inserting \eqref{lemmaA3} into \eqref{lemmaA2}, we can get
\begin{align}\label{lemmaA4}
&\|\Box_{m,n}((A(x,D)-A(m,n))u)\|_r\nonumber\\
&\lesssim\sum_{h,k,l\in\mathbb{Z}^d}\langle m-h\rangle^{-{N_1}}\langle n-l\rangle^{-{N_2}}\langle h-k\rangle^{-{N_3}}\big(\langle l\rangle^{M-\varepsilon}+\langle n\rangle^{M-\varepsilon}\big)\|\overline{\Box}_{k,l}u\|_r.
\end{align}
We divide \eqref{lemmaA4} into $\langle l\rangle\leq 10 \langle n\rangle$ and $\langle l\rangle > 10\langle n\rangle$ two parts,
\begin{align}\label{lemmaA5}
&\|\Box_{m,n}((A(x,D)-A(m,n))u)\|_r\nonumber\\
&\lesssim\sum_{h,k,l\in\mathbb{Z}^d}\langle m-h\rangle^{-{N_1}}\langle n-l\rangle^{-{N_2}}\langle h-k\rangle^{-{N_3}}\langle n\rangle^{M-\varepsilon}\|\overline{\Box}_{k,l}u\|_r\chi_{\{\langle l\rangle\leq 10 \langle n\rangle\}}\nonumber\\
&\quad +\sum_{h,k,l\in\mathbb{Z}^d}\langle m-h\rangle^{-{N_1}}\langle n-l\rangle^{-{N_2}+M-\varepsilon}\langle h-k\rangle^{-{N_3}}\|\overline{\Box}_{k,l}u\|_r\chi_{\{\langle l\rangle> 10 \langle n\rangle\}}.
\end{align}
 From Young's inequality, we have
\begin{align}\label{lemmaA6a}
&\|\mathscr{A}_{m,n}((A(x,D)-A(m,n))u)\|_{L^\gamma(0,T;L^r)}\nonumber\\
&\lesssim \|e^{-tA(m,n)}\|_{L^1(0,T)} \|\Box_{m,n}((A(x,D)-A(m,n))u)\|_{L^\gamma(0,T;L^r)},
\end{align}
and
\begin{align}\label{lemmaA6b}
&\|\mathscr{A}_{m,n}((A(x,D)-A(m,n))u)\|_{ L^\infty(0,T;L^r)}\nonumber\\
&\lesssim \|e^{-tA(m,n)}\|_{L^{\gamma'}(0,T)} \|\Box_{m,n}((A(x,D)-A(m,n))u)\|_{L^\gamma(0,T;L^r)}.
\end{align}
For $\langle l\rangle\leq 10 \langle n\rangle$, inserting \eqref{lemmaA5} into \eqref{lemmaA6a} and \eqref{lemmaA6b}, we know from Young's inequality that
\begin{align}\label{lemmaA9a}
&\Big\|\|\mathscr{A}_{m,n}((A(x,D)-A(m,n))u)\|_{L^\gamma(0,T;L^r)}\Big\|_{\ell^p_m}\nonumber\\
&\lesssim  \big(\langle n\rangle^{M-\varepsilon}\|e^{-tA(m,n)}\|_{L^1(0,T)}\big)\nonumber\\
&\ \ \ \ \times \sum_{l\in\mathbb{Z}^d}\langle n-l\rangle^{-{N_2}}\bigg\|\sum_{h,k\in\mathbb{Z}^d}\langle m-h\rangle^{-{N_1}}\langle h-k\rangle^{-{N_3}}\|\overline{\Box}_{k,l}u\|_{L^\gamma(0,T;L^r)}\bigg\|_{\ell^p_m}\chi_{\{\langle l\rangle\leq 10 \langle n\rangle\}}\nonumber\\
&\lesssim \big(\langle n\rangle^{M-\varepsilon}\|e^{-tA(m,n)}\|_{L^1(0,T)}\big)\cdot\sum_{l\in\mathbb{Z}^d}\langle n-l\rangle^{-{N_2}}\big\|\|\overline{\Box}_{m,l}u\|_{L^\gamma(0,T;L^r)}\big\|_{\ell^p_m}\chi_{\{\langle l\rangle\leq 10 \langle n\rangle\}},
\end{align}
and
\begin{align}\label{lemmaA9b}
&\Big\|\|\mathscr{A}_{m,n}((A(x,D)-A(m,n))u)\|_{L^\infty(0,T;L^r)}\Big\|_{\ell^p_m}\nonumber\\
&\lesssim \big(\langle n\rangle^{M-\varepsilon}\|e^{-tA(m,n)}\|_{L^{\gamma'}(0,T)}\big) \sum_{l\in\mathbb{Z}^d}\langle n-l\rangle^{-{N_2}}\big\|\|\overline{\Box}_{m,l}u\|_{L^\gamma(0,T;L^r)}\big\|_{\ell^p_m}\chi_{\{\langle l\rangle\leq 10 \langle n\rangle\}}.
\end{align}
Therefore,
\begin{align}\label{lemmaA7}
&\|\mathscr{A}_{m,n}((A(x,D)-A(m,n))u)\|_{\mathscr{L}^\gamma(0,T;\ell^q_{s_\gamma} \ell^p(L^r))\cap \mathscr{L}^\infty(0,T;\ell^q_{s_0} \ell^p(L^r))}\nonumber\\
&\lesssim\int^T_0 \sup_{n\in\mathbb{Z}^d}\frac{\langle n\rangle^{M-\varepsilon}}{e^{t|n|^M}}dt+\bigg(\int^T_0 \Big(\sup_{n\in\mathbb{Z}^d}\frac{\langle n\rangle^{M-M/\gamma-\varepsilon}}{e^{t|n|^M}}\Big)^{\gamma'}dt\bigg)^{1/\gamma'}\nonumber\\
&\ \ \ \times \Big\|\langle n\rangle^{s_\gamma}\sum_{l\in\mathbb{Z}^d}\langle n-l\rangle^{-{N_2}}\big\|\|\overline{\Box}_{m,l}u\|_{L^\gamma(0,T;L^r)}\big\|_{\ell^p_m}\chi_{\{\langle l\rangle\leq 10 \langle n\rangle\}}\Big\|_{\ell^q_n}.
\end{align}
Noticing $s_\gamma\geq 0$, taking $N_2>s_\gamma+d$, and using Young's inequality, we obtain
\begin{align}\label{lemmaA8}
&\|\mathscr{A}_{m,n}((A(x,D)-A(m,n))u)\|_{\mathscr{L}^\gamma(0,T;\ell^q_{s_\gamma} \ell^p(L^r))\cap \mathscr{L}^\infty(0,T;\ell^q_{s_0} \ell^p(L^r))}\nonumber\\
&\lesssim \big(T+T^{\varepsilon/M}+T^{1/\gamma'}\big) \Big\|\sum_{l\in\mathbb{Z}^d}\langle n-l\rangle^{-{N_2}+{s_\gamma}}\langle l\rangle^{s_\gamma}\big\|\|\overline{\Box}_{m,l}u\|_{L^\gamma(0,T;L^r)}\big\|_{\ell^p_m}\Big\|_{\ell^q_n}\nonumber\\
&\lesssim \big(T+T^{\varepsilon/M}+T^{1/\gamma'}\big)\|u\|_{\mathscr{L}^\gamma(0,T;X^{s_\gamma}_{r,p,q})}.
\end{align}
For $\langle l\rangle> 10 \langle n\rangle$, we know from Young's inequality that
\begin{align}
&\Big\|\|\mathscr{A}_{m,n}((A(x,D)-A(m,n))u)\|_{L^\gamma(0,T;L^r)\cap L^\infty(0,T;L^r) }\Big\|_{\ell^p_m}\nonumber\\
&\lesssim \big(T+T^{1/\gamma'}\big) \sum_{l\in\mathbb{Z}^d}\langle n-l\rangle^{-{N_2}+M-\varepsilon}\bigg\|\sum_{h,k\in\mathbb{Z}^d}\langle m-h\rangle^{-{N_1}}\langle h-k\rangle^{-{N_3}}\|\overline{\Box}_{k,l}u\|_{L^\gamma(0,T;L^r)}\bigg\|_{\ell^p_m}\nonumber\\
&\lesssim \big(T+T^{1/\gamma'}\big)\sum_{l\in\mathbb{Z}^d}\langle n-l\rangle^{-{N_2}+M-\varepsilon}\big\|\|\overline{\Box}_{m,l}u\|_{L^\gamma(0,T;L^r)}\big\|_{\ell^p_m}.
\end{align}
Then using the similar method with \eqref{lemmaA7}-\eqref{lemmaA8}, just taking $N_2>M+s_\gamma+d$, we can get
\begin{align}\label{lemmaA10}
&\|\mathscr{A}_{m,n}((A(x,D)-A(m,n))u)\|_{\mathscr{L}^\gamma(0,T;\ell^q_{s_\gamma} \ell^p(L^r))\cap \mathscr{L}^\infty(0,T;\ell^q_{s_0} \ell^p(L^r))}\nonumber\\
&\lesssim \big(T+T^{1/\gamma'}\big)\|u\|_{\mathscr{L}^\gamma(0,T;X^{s_\gamma}_{r,p,q})}.
\end{align}
Now the proof is completed. $\hfill\Box$

\begin{lem}\label{lemmaC} Let  $1\leq p,q,r\leq\infty$. Assume that $A(x,\xi)$ satisfies (A1) and (A2), and $A(0,0)>0$.  Then we have
\begin{gather}
\|\mathscr{A}_{m,n}((A(x,D)-A(m,n))u)\|_{\mathscr{L}^2(0,\infty;\ell^q_{s_2} \ell^p(L^r))\cap\mathscr{L}^\infty(0,\infty;\ell^q_{s_0} \ell^p(L^r))}\lesssim \|u\|_{\mathscr{L}^2(0,\infty;X^{s_2-\varepsilon}_{r,p,q})}.\label{lemmaCa}
\end{gather}
\end{lem}
{\bf Proof.} One just need to note that
\begin{align*}
&\|\mathscr{A}_{m,n}((A(x,D)-A(m,n))u)\|_{L^2_tL_x^r}\nonumber\\
&\lesssim \int^\infty_0 e^{-tA(m,n)}dt\cdot\|\Box_{m,n}((A(x,D)-A(m,n))u)\|_{L^2_tL_x^r}\nonumber\\
&\lesssim \langle n\rangle^{-M}\|\Box_{m,n}((A(x,D)-A(m,n))u)\|_{L^2_tL_x^r},
\end{align*}
and
\begin{align*}
 \|\mathscr{A}_{m,n}((A(x,D)-A(m,n))u)\|_{L^\infty_tL_x^r}
 \lesssim \langle n\rangle^{-M/2}\|\Box_{m,n}((A(x,D)-A(m,n))u)\|_{L^2_tL_x^r}.
\end{align*}
Taking the similar proceeding with \eqref{lemmaA9a}-\eqref{lemmaA10}, we can get the conclusion \eqref{lemmaCa}. $\hfill\Box$\\

\section{Appendix}
In this Appendix we show that the solution of the system \eqref{localization integral} is really solving \eqref{PDE} in some function spaces. Denote
\begin{align}\label{spaceapp0}
  \|f\|_{X_{r,p,q}^{s_1,s_2}}= \left\|\big\| \{\|\langle m\rangle^{s_1} \langle n\rangle^{s_2}\Box_{m,n} f \|_{r}\}  \big\|_{\ell^p_{m\in \mathbb{Z}^d}} \right\|_{\ell^q_{n\in \mathbb{Z}^d}}.
\end{align}
Comparing $X_{r,p,q}^{s_1,s_2}$ with $X_{r,p,q}^{s_2}$, we see that $X_{r,p,q}^{s_1,s_2}$ contains another weight $\langle m\rangle^{s_1}$ with respect to physical spaces. We further denote
\begin{align}\label{spaceapp10}
  \|f\|_{\mathscr{L}^{\gamma}(0,T; X_{r,p,q}^{s_1,s_2})}= \left\| \big\| \{\|\langle m\rangle^{s_1} \langle n\rangle^{s_2}\Box_{m,n} f \|_{L^\gamma(0,T;L^r)} \} \big\|_{\ell^p_{m\in \mathbb{Z}^d}} \right\|_{\ell^q_{n\in \mathbb{Z}^d}}.
\end{align}
For the sake of convenience, we also use the notation
\begin{align}\label{spaceapp100}
  \|\{g_{m,n}\}\|_{\mathscr{L}^{\gamma}(0,T; \ell^q_{s_2} \ell^p_{s_1} (L^r))}= \left\| \big\| \{\|\langle m\rangle^{s_1} \langle n\rangle^{s_2} g_{m,n}   \|_{L^\gamma(0,T;L^r)} \} \big\|_{\ell^p_{m\in \mathbb{Z}^d}} \right\|_{\ell^q_{n\in \mathbb{Z}^d}}.
\end{align}

\begin{prop}
Let $s_c >\kappa+ d/q' $,   $\gamma\gg \gamma(K) \vee \sigma_1 \vee \sigma_2$ is defined in Theorem \ref{thma}.  Then the solution obtained in Theorem \ref{thma} satisfies \eqref{PDE} in the sense of
$$
\partial_t u, \ A(x,D) u , \ F(u) \in  \mathscr{L}^\gamma(0,T; X^{-\sigma_1, s_c   - \sigma_2}_{r,p,q})
$$
\end{prop}
{\bf Proof.} Taking the derivative on time variable to  \eqref{Estimate22}, we have
\begin{align}\label{localizationdiff}
&\partial_t \Box_{m,n} u^{(\mu+1)} = -    \Box_{m,n}(A(x,D) u^{(\mu+1)}) + \Box_{m,n}F(u^{(\mu)}),
  \  \ \ \    {\rm for \ all\ } m,n \in \mathbb{Z}^d.
\end{align}
Let us observe that
\begin{align} \label{decom}
 \Box_{m,n}(A(x,D) u^{(\mu+1)})=  A(m,n) \Box_{m,n} u^{(\mu+1)}  + \Box_{m,n} \mathscr{F}^{-1}(A(x,\xi)-A(m,n))\widehat{u}^{(\mu+1)}.
\end{align}
It is easy to see that
\begin{align} \label{app1}
\left\|\{A(m,n) \Box_{m,n} (u^{(\mu+1)}-u)\} \right\|_{\mathscr{L}^{\gamma}(0,T; \ell^q_{s_c -\sigma_2} \ell^p_{-\sigma_1} (L^r))}  \lesssim \|u^{(\mu+1)} -u\|_{\mathscr{L}^{\gamma}(0,T; X_{r,p,q}^{s_c})} \to 0.
\end{align}
Moreover, by Lemma \ref{lemma4} we have
\begin{align} \label{app2}
\left\| \{\Box_{m,n}((a(x) -a(m)) (u^{(\mu+1)}-u))\} \right\|& _{\mathscr{L}^{\gamma}(0,T; \ell^q_{s_c -\sigma_2} \ell^p_{-\sigma_1} (L^r))}  \nonumber\\
&   \lesssim \|u^{(\mu+1)} -u\|_{\mathscr{L}^{\gamma}(0,T; X_{r,p,q}^{s_c})} \to 0.
\end{align}
It follows from Lemma \ref{lemma5} that
\begin{align} \label{app3}
\| \{\Box_{m,n}((b(D) -b(n)) (u^{(\mu+1)}-u))\} \|& _{\mathscr{L}^{\gamma}(0,T; \ell^q_{s_c -\sigma_2} \ell^p_{-\sigma_1} (L^r))}  \nonumber\\
&  \lesssim \|u^{(\mu+1)} -u\|_{\mathscr{L}^{\gamma}(0,T; X_{r,p,q}^{s_c})} \to 0.
\end{align}
In view of the algebra property
\begin{align} \label{app4}
& \left\| \{\Box_{m,n}(F(u^{(\mu)}) - F(u))\} \right\|_{\mathscr{L}^{\gamma}(0,T; \ell^q_{s_c -\sigma_2} \ell^p_{-\sigma_1} (L^r))}  \nonumber\\
& \quad \lesssim \omega \left(\|u^{(\mu)}\|_{\mathscr{L}^{\infty}(0,T; X_{r,p,q}^{s_c})}, \ \|u\|_{\mathscr{L}^{\infty}(0,T; X_{r,p,q}^{s_c})} \right) \|u^{(\mu)} -u\|_{\mathscr{L}^{\gamma}(0,T; X_{r,p,q}^{s_c})} \to 0.
\end{align}
Collecting the estimates in \eqref{app1}--\eqref{app4}, we see that
$$
 -A(x,D) u^{(\mu+1)}- +  F(u^{(\mu)})   \to  - A(x,D) u  +  F(u )
$$
in $\mathscr{L}^\gamma(0,T; X^{-\sigma_1, s_c   - \sigma_2}_{r,p,q})$.  By \eqref{localizationdiff}, we see that
$  \partial_t u^{(\mu+1)} $ converges to $ A(x,D) u  +  F(u )$ in $\mathscr{L}^\gamma(0,T; X^{-\sigma_1, s_c   - \sigma_2}_{r,p,q})$. Using the embedding
$$
\mathscr{L}^\gamma(0,T; X^{-\sigma_1, s_c  - \sigma_2}_{r,p,q}) \subset \mathscr{D}'(0,T; \mathscr{S}'(\mathbb{R}^d)),
$$
we see that $  \partial_t u^{(\mu+1)} \to \partial_t u $. So, $u$ satisfies \eqref{PDE} in  $\mathscr{L}^\gamma(0,T; X^{-\sigma_1, s_c   - \sigma_2}_{r,p,q})$. $\hfill\Box$\\

If the initial data have more regularity, we can show that

\begin{prop}
Let $s_c > \kappa+ d/q' $, $s_\gamma =  s_c+ \sigma_2/\gamma$, $\gamma$ is defined in Theorem \ref{thma}. Let $u_0\in M^{s_c +\sigma_2}_{p,q} \cap X^{\sigma_1,s_c}_{r,p,q}$. Then the solution obtained in Theorem \ref{thma} satisfies \eqref{PDE} in the sense of
$$
\partial_t u, \ A(x,D) u , \ F(u) \in  \mathscr{L}^\gamma(0,T; X^{s_c + \sigma_2/\gamma}_{r,p,q}) \cap \mathscr{L}^\infty (0,T; X^{s_c}_{r,p,q}).
$$

\end{prop}
{\bf Proof.} First, we estimate $\partial_t u$. According to the  system  \eqref{limit22},
\begin{align}\label{limit101}
 \Box_{m,n}  \partial_t u (t) = & - A(m,n) e^{-t A(m,n)}\Box_{m,n} u_0 -  e^{-t A(m,n)}\Box_{m,n}((a(x)-a(m))   u_0 ) \nonumber\\
 &  - e^{-t A(m,n)}\Box_{m,n} ((b(D)-b(n)) u_0 ) +  e^{-t A(m,n)}\Box_{m,n} F(u_0) \nonumber\\
  &-   \mathscr{A}_{m,n}((a(x)-a(m)) \partial_t u ) -   \mathscr{A}_{m,n}((b(D)-b(n)) \partial_t u )+   \mathscr{A}_{m,n} \partial_t (F(u )) \nonumber\\
  :=& I_{m,n}+...+VII_{m,n}.
\end{align}
Using \eqref{semigroupest1} and $|A(m,n)| \lesssim \langle m\rangle^{\sigma_1} + \langle n\rangle^{\sigma_2}$, we see that
\begin{align}\label{limit102}
 \|\{ I_{m,n} \}\|_{\mathscr{L}^\gamma(0,T;\ell^q_{s_\gamma} \ell^p(L^r))}  \lesssim \|u_0\|_{M^{s_c+\sigma_2}_{p,q}} +  \|u_0\|_{X^{\sigma_1, s_c}_{r,p,q}}.
\end{align}
By Lemma \ref{lemma4} and \eqref{semigroupest1}, we have
\begin{align}\label{limit10a}
&\|e^{-tA(m,n)}\Box_{m,n}((a(x)-a(m))u_0)\|_{L^{\gamma}(0,T; L^r)}\nonumber\\
&\lesssim  \langle n \rangle^{ -\sigma_2/\gamma } \sum_{n_1,n_2,l\in\mathbb{Z}^d}\langle m-l\rangle^{-N+\sigma_1}\langle l\rangle^{ \sigma_1 }\langle n_1\rangle^{-K}\|\widetilde{\Box}_{l,n_2}u_0\|_r \chi_{|n-n_1-n_2|_\infty \leq 3}.
\end{align}
Using Young's inequality, we have
\begin{align}\label{limit11a}
 \|\{ II_{m,n} \}\|_{\mathscr{L}^\gamma(0,T;\ell^q_{s_\gamma} \ell^p(L^r))}  \lesssim  \|u_0\|_{X^{\sigma_1, s_c}_{r,p,q}}.
\end{align}
Applying Lemma \ref{lemma5} and  \eqref{semigroupest1},
\begin{align}\label{limit12a}
 \|e^{-tA(m,n)} & \Box_{m,n}((b(D)-b(n))u_0)\|_{L^{\gamma}(0,T; L^r)}  \nonumber\\
  & \lesssim  \langle n \rangle^{ -\sigma_2/\gamma +\sigma_2-1}\sum_{|l_1|_\infty\leq 1}\sum_{k\in\mathbb{Z}^d}\langle m-k\rangle^{-(d+1-\vartheta)} \|\Box_{k,n+l_1}u_0 \|_r.
\end{align}
Applying Young's inequality, we have
\begin{align}\label{limit13a}
 \|\{ III_{m,n} \}\|_{\mathscr{L}^\gamma(0,T;\ell^q_{s_\gamma} \ell^p(L^r))}  \lesssim  \|u_0\|_{M^{\sigma_2 + s_c}_{p,q}}.
\end{align}
Since $M^{s_0}_{p,q}$ is an algebra, by \eqref{semigroupest1} we have\footnote{For the general nonlinearity, we have similar result.}
\begin{align}\label{limit14a}
 \|\{ IV_{m,n} \}\|_{\mathscr{L}^\gamma(0,T;\ell^q_{s_\gamma} \ell^p(L^r))}    \lesssim  \|F(u_0)\|_{M^{ s_c}_{p,q}}
 \lesssim  \|u_0\|^K_{M^{ s_c+\sigma_2}_{p,q}} + \|u_0\|^L_{M^{ s_c+\sigma_2}_{p,q}}.
\end{align}
The estimates of $V_{m,n}-VII_{m,n}$ are similar to those in Lemmas \ref{locallemma2}, \ref{locallemma3} and \ref{locallemma4}:
\begin{align}
 \|\{ V_{m,n} \}\|_{\mathscr{L}^\gamma(0,T;\ell^q_{s_\gamma} \ell^p(L^r))}   &  \lesssim   T^\delta \|\partial_t u\|_{\mathscr{L}^\gamma(0,T; X^{s_\gamma}_{r,p,q})},   \label{limit15a}\\
  \|\{ VI_{m,n} \}\|_{\mathscr{L}^\gamma(0,T;\ell^q_{s_\gamma} \ell^p(L^r))}   &  \lesssim   T^\delta \|\partial_t u\|_{\mathscr{L}^\gamma(0,T; X^{s_\gamma}_{r,p,q})},   \label{limit16a}\\
   \|\{ VII_{m,n} \}\|_{\mathscr{L}^\gamma(0,T;\ell^q_{s_\gamma} \ell^p(L^r))}   &  \lesssim   T^\delta (\| u\|^{K-1}_{\mathscr{L}^\gamma(0,T; X^{  s_\gamma}_{r,p,q})} + \| u\|^{L-1}_{\mathscr{L}^\gamma(0,T; X^{ s_\gamma}_{r,p,q})}) \nonumber\\
   & \ \ \ \ \times\|\partial_t u\|_{\mathscr{L}^\gamma(0,T; X^{  s_\gamma}_{r,p,q})}.   \label{limit17a}
\end{align}
Collecting the estimates of $I_{m,n}-VII_{m,n}$  and using the same way as the argument in proving the local well-posedness, we have for some $0<T<1$,
$$
\partial_t u \in  \mathscr{L}^\gamma(0,T; X^{s_\gamma}_{r,p,q}) \cap \mathscr{L}^\infty (0,T; X^{s_c}_{r,p,q}).
$$
By Lemmas \ref{lemma4} and \ref{lemma5}, we have
$$
A(x,D) u \in  \mathscr{L}^\gamma(0,T; X^{s_\gamma}_{r,p,q}) \cap \mathscr{L}^\infty (0,T; X^{s_c}_{r,p,q}).
$$
Similarly as in the above,
$$
F( u) \in  \mathscr{L}^\gamma(0,T; X^{s_\gamma}_{r,p,q}) \cap \mathscr{L}^\infty (0,T; X^{s_c}_{r,p,q}).
$$
So, we have the result, as desired. $\hfill\Box$\\

\noindent{\bf Acknowledgment.} The second named author is grateful to Professor Hans Feichtinger for his valuable discussions on modulation spaces with respected to the paper \cite{DoFeGr06}. The first and second named authors are
supported in part by the National Science Foundation of China, grants   11271023. The third named author is partially supported by NSFC 11626041 and Beijing Science Foundation 1162004.  The fourth named author has been partially supported by the Natural Sciences and Engineering Research Council of Canada under Discovery Grant 0008562.

\footnotesize

\end{document}